\documentclass[3p, number,sort&compress,times,fleqn]{elsarticle}

\usepackage{amsmath}
\usepackage{amssymb}
\usepackage{mathtools}
\usepackage[colorlinks = true,
linkcolor = blue,
urlcolor  = blue,
citecolor = blue,
anchorcolor = blue]{hyperref}
\usepackage[utf8]{inputenc}
\usepackage{subfig}
\usepackage{caption}
\usepackage{mathdots}
\usepackage{lineno}

\usepackage{csml}
\usepackage{tikz}
\usepackage{quantikz}
\usepackage{braket}
\usepackage{ulem}
\usepackage{relsize}

\begin{document}
	\begin{frontmatter}
	
	\title{A quantum spectral method for non-periodic boundary value problems}
	\author[1]{Eky Febrianto}
	\ead{eky.febrianto@glasgow.ac.uk}
	\author[2]{Yiren Wang}
	\ead{yw680@cam.ac.uk}
	\author[2]{Burigede Liu}
	\ead{bl377@cam.ac.uk}
	\author[3]{Michael Ortiz}
	\ead{ortiz@aero.caltech.edu}	
	\author[2]{Fehmi Cirak \corref{cor1}}
	\ead{fc286@cam.ac.uk}
		
	\cortext[cor1]{Corresponding author}
	
	\address[1]{Glasgow Computational Engineering Centre, James Watt School of Engineering, University of Glasgow, Glasgow, G12 8QQ, UK}
	\address[2]{Department of Engineering, University of Cambridge, Cambridge, CB2 1PZ, UK}
	\address[3]{Division of Engineering and Applied Science, California Institute of Technology, Pasadena, CA 91125, USA}
	
	\begin{abstract}
	Quantum computing offers the potential to solve computational mechanics problems with polylogarithmic complexity $\mathcal{O}((\log N)^c)$. We present a quantum spectral method for non-periodic boundary value problems with arbitrary Dirichlet boundary conditions. Our method extends the recently proposed approach by Liu {\rm et al.} (2024), in which periodic problems are discretised using truncated Fourier series. In such methods, the discretisation of boundary value problems with constant coefficients leads to a set of algebraic equations in Fourier space. We implement the associated diagonal solution operator by first approximating it with a polynomial and then quantum encoding it. The mapping between the physical and Fourier spaces is accomplished using the quantum Fourier transform (QFT). To impose zero Dirichlet boundary conditions, we double the domain size and reflect all physical fields antisymmetrically. The reflection defines the quantum sine transform (QST) by pre- and post-multiplying with the QFT. For non-zero Dirichlet boundary conditions, we decompose the solution into a boundary-conforming and a homogeneous part. The homogeneous part is determined by solving a problem with a suitably modified forcing vector. We illustrate the proposed approach with a Dirichlet-Poisson problem and demonstrate its generality by applying it to a fractional stochastic PDE for modelling spatial random fields. We discuss the implementation of the quantum circuits and provide numerical evidence confirming polylogarithmic complexity. A Qiskit implementation of all circuits is available on GitHub.	
\end{abstract}
	
	
\end{frontmatter}

\section{Introduction}
\label{sec:introduction}
%
The solution of boundary value problems from computational mechanics is one of the most promising application areas of future quantum computing~\cite{liu2024towards}. The key appeal of quantum computing is that problems can be solved in $\mathcal {O}(\log^c N)$ time, instead of $\mathcal{O}( N^c)$ in classical computing, where $N$ is the problem size and $c>1$ a constant. In quantum computing, the $N$-dimensional forcing and solution vectors of a discretised boundary value problem are encoded as the complex-valued state vector of a system with~$n=\log N$ quantum particles, or qubits ($\log$ denotes~$\log_2$). The state vector of each qubit is of dimension two. The solution of the problem is formulated as the mapping of the forcing state vector to the solution state vector using a suitably constructed length and angle preserving unitary matrix,  encoding the discretised solution operator. Quantum programming concerns the assembly of the problem-dependent unitary matrix from simpler elementary unitary matrices provided by the quantum hardware, considering that the set of unitary matrices is closed under multiplication ~\cite{deutsch1989quantum,barenco1995elementary}. The elementary unitary matrices, or gates, operate on a single qubit or at most two qubits at a time. Informally, quantum speed-up is achieved by designing algorithms that have a polynomial number of one or two-qubit gates operating on the~$n$ qubits while the vectors and matrices operated on live in the Hilbert space~$\mathbb C^N$ of dimension~$N=2^n$. The  exponentially compressed quantum representation makes it very costly to address and modify selected coefficients of a vector or a matrix, or to implement constructs like the if-then-else commonly used in classical computing. 

It is evident that quantum computing is fundamentally different from classical computing. Hence, the present finite element approaches, primarily based on unstructured meshes and locally supported basis functions, may not be well-suited for quantum computing. For unstructured problems, preparing a desired quantum state on a quantum computer usually incurs a cost of~$\mathcal {O}(N)$, invalidating any potential quantum advantage in subsequent processing of the data. For instance, in the case of state preparation, which is an extensively studied problem, encoding a classical $N$ dimensional vector has complexity~$\mathcal{O}(N)$ when no ancillary qubits are used, see e.g.~\cite{mottonen2004transformation,sun2023asymptotically}. The same applies to the quantum encoding of a matrix, which can be thought of as~$N$ vectors of dimension $N$ or less, when the matrix is sparse. The computational cost reduces to~$\mathcal{O}(\log N)$ when the problem can be formulated in a way that allows vectors or matrices to be approximated by a function~\cite{rosenkranz2025quantum}. Setting aside the cost of state preparation, there are several algorithms, most prominently, Quantum Singular Value Transformation (QSVT)~\cite{martyn2021grand} and HHL~\cite{harrow2009quantum}, for solving linear systems of equations with polylogarithmic complexity. The original HHL algorithm has a complexity~$O(\kappa^2/\epsilon \log N)$ and QSVT has a complexity~$O(\kappa \log (\kappa/\epsilon))$, where $\kappa$ is the condition number and $\epsilon$ the accuracy of the solution~\cite{morales2024quantum}. These and other complexity estimates from quantum computing assume that all vectors and matrices can be accessed via an oracle, i.e. a subroutine. The complexity of the oracle implementation is not taken into account. The quantum encoding of (sparse) matrices, and the oracles, is discussed in more recent works~\cite{camps2022fable,camps2024explicit,lapworth2024evaluation,sunderhauf2024block}. Regardless of the cost of quantum encoding, the available linear system solvers, including QSVT and HHL, either require a large number of quantum operations or demand a large number of ancilla qubits. This is exacerbated by the fact that the qubits must be fault-tolerant, which typically requires redundancy and replaces one (logical) qubit with several physical qubits. Consequently, the resulting hardware requirements for practically relevant problem sizes exceed what is expected to be available in the near term. 

We conjecture that the simple structure and straightforward implementation of classical spectral methods on uniform grids make them especially appealing for quantum computing; see~\cite{trefethen2000spectral, canuto2007spectral} for an introduction to spectral methods.   The commonly used basis functions in spectral methods are the trigonometric Fourier basis functions for periodic problems and Chebyshev basis functions for non-periodic problems. Both kinds of basis functions have appealing orthogonality properties and are global in the sense that they are supported over the entire problem domain. The discretisation of spectral methods leads to well-structured vectors and matrices with repetitive entries, which are easier to quantum encode than their fully unstructured equivalents from other discretisation methods. In particular, the Fourier-based approaches lead in the Fourier space to diagonal matrices or block-diagonal matrices when the coefficients of the differential equation are constant. The diagonality is independent of the dimension of the problem domain and the order of the differential operators, whether integer or fractional.  Indeed, as recently demonstrated in Liu et al.~\cite{liu2024towards}, for FFT-based homogenisation~\cite{moulinec1998numerical,schneider2021review,gierden2022review}, it is possible to solve specific problems from solid mechanics using a Fourier-based spectral approach with complexity~$\mathcal{O}(\log^c N$).  Recently, the feasibility of spectral methods on currently available quantum hardware has also been demonstrated in~\cite{wright2024noisy} for a prototypical one-dimensional wave propagation problem.   There are numerous works in quantum computing literature on spectral-like methods based on  Fourier basis functions~\cite{childs2021high}, Chebyshev basis functions~\cite{wu2025quantum}, and finite differences on structured grids~\cite{cao2013quantum,lubasch2025quantum}. The periodic/circulant system matrices in finite differences can be diagonalised using a discrete Fourier transform (DFT), or its equivalent, the quantum Fourier transform (QFT), leading to a spectral-like approach~\cite[Ch. 5]{strang1986introduction}. The mentioned works mainly differ in terms of the problems investigated and in how the resulting systems of equations with diagonal or block-diagonal matrices are solved. 

In this work, we introduce a novel quantum spectral solver for boundary value problems with arbitrary Dirichlet boundary conditions on axis-aligned~$d$-dimensional hyperrectangles.  We illustrate the basic approach with a Dirichlet-Poisson problem and demonstrate its generality by applying it to a fractional stochastic PDE used in modelling spatial random fields~\cite{lindgren2011explicit,lindgren2022spde,koh2023stochastic,ben-yelun2024struct}. Our method builds on the foundational work by Liu et al.~\cite{liu2024towards} on solving periodic Poisson problems in homogenisation. In~\cite{liu2024towards}, the periodic fields are discretised using a truncated Fourier basis, and the corresponding coefficients are determined by solving independent scalar equations in the Fourier space. The mapping between physical and Fourier spaces is carried out using the quantum Fourier transform (QFT), which has complexity~$\mathcal{O}(\log^2 N)$, in contrast to the classical DFT with complexity~$\mathcal{O}(N \log N)$. The mode-wise (frequency-wise) scalar division in Fourier space is a non-length-preserving operation that can only be implemented approximately on a quantum computer. We implement scalar division via polynomial approximation of the components of the diagonal solution matrix. The specific polynomial encoding technique we employ is the multivariate extension of the univariate scheme introduced in~\cite{woerner2019quantum}. We extend the approach of Liu et al.~\cite{liu2024towards} to non-periodic problems by doubling the physical domain size and antisymmetrically reflecting the source term. The new source and solution state vectors require only~$d$ additional qubits compared to those for the original domain, where~$d$ is the spatial dimension. The application of the QFT to the antisymmetrically extended source state vector cancels all cosine terms in the Fourier series, leaving only boundary conforming sine terms.  Formally, a unitary reflection matrix is introduced to define the quantum sine transform (QST), the analogue of the discrete sine transform (DST), by pre- and post-multiplying the QFT matrix with the reflection matrix~\cite{wickerhauser1996adapted,strang1999discrete}. Our quantum implementation of the reflection matrix follows and improves upon~\cite{klappenecker2001discrete}. We consider different implementations of its components, specifically the forward shift operator~\cite{camps2024explicit,wang2025comprehensive}, to reduce gate counts. In summary, replacing the QFT in the periodic PDE solver with DST yields a solver for homogeneous Dirichlet problems. To consider inhomogeneous boundary conditions, we use the superposition principle and decompose the solution into a boundary conforming part and a homogeneous part. The boundary conforming part is a sufficiently smooth function that can be chosen freely, while the homogeneous part is obtained by solving the homogeneous PDE with a modified source vector. A Qiskit implementation of the proposed quantum spectral solver is provided on GitHub~\cite{qspectral2026}.

%
\section{Classical Fourier spectral method \label{sec:classical-fsm}}
%
In this section, we provide a brief review of the classical Fourier spectral method for solving partial differential equations and the treatment of homogeneous and inhomogeneous Dirichlet boundary conditions. 
For clarity, we consider a one-dimensional Poisson-Dirichlet problem with homogeneous boundary conditions, followed by its generalisation to other differential equations. In particular, we study a fractional differential equation motivated by the description of random fields. Finally, we introduce the treatment of inhomogeneous boundary conditions. 
%
\subsection{Poisson-Dirichlet problem \label{sec:model-problem}}
\begin{figure}[]
\centering
\begin{tikzpicture}
    \draw[thick] (0,0) -- (12,0);
    \foreach \x in {0,1,2,3,5,6,7} {
        \draw[thin] (1.5*\x,0.1) -- (1.5*\x,-0.1); 
        \node[below] at (1.5*\x,-0.1) {$\x$}; 
    }
    
     \draw[thin] (1.5*4,0.1) -- (1.5*4,-0.1); 
     \node[below] at (1.5*4,-0.1) {$N/2=4$}; 
      \node[below] at (1.5*8,-0.1) {$0$}; 
    
    \draw[domain=0:6,smooth,variable=\x,csmlBrightBlue,line width=0.25mm] 
        plot ({\x}, {(0.3333*\x +3/4)}); 
    \draw[domain=6:12,smooth,variable=\x,csmlBrightBlue,line width=0.25mm] 
        plot ({\x}, {(0.33333*\x-4)-3/4 }); 
		\fill[black!60](0.0,0.0) circle[radius=0.1];
		\fill[csmlBrightBlue] (1.5,0.5*1+3/4) circle[radius=0.1];
		\fill[csmlBrightBlue] (3,0.5*2+3/4) circle[radius=0.1];
		\fill[csmlBrightBlue] (4.5,0.5*3+3/4) circle[radius=0.1];
		\fill[black!60] (6.0,0.0) circle[radius=0.1];
		\fill[csmlBrightBlue] (7.5,-0.5*3-3/4) circle[radius=0.1];
		\fill[csmlBrightBlue] (9,-0.5*2-3/4) circle[radius=0.1];
		\fill[csmlBrightBlue] (10.5,-0.5*1-3/4) circle[radius=0.1];
		\fill[black!60]  (1.5*8, 0) circle[radius=0.1];
    \node[above] at (0,0.15) {$x=0$};
    \node[above] at (6,0.15) {$x=L$};
    \node[above] at (12,0.15) {$x=2 L$};
\end{tikzpicture}
\caption{The extended domain~$\Omega_{\text E}=(0, \, 2L)$ and its discretisation with~$N = 8$ cells and equidistant grid points. The original problem domain is~$\Omega=(0, \, L)$.  A sample reflected source term~$f_\text{E}(x)$ is depicted in blue, and the components of the respective force vector~$\vec f$ are depicted as dots. The extended domain problem is~$2L$ periodic.\label{fig:extended_domain}}
\end{figure}
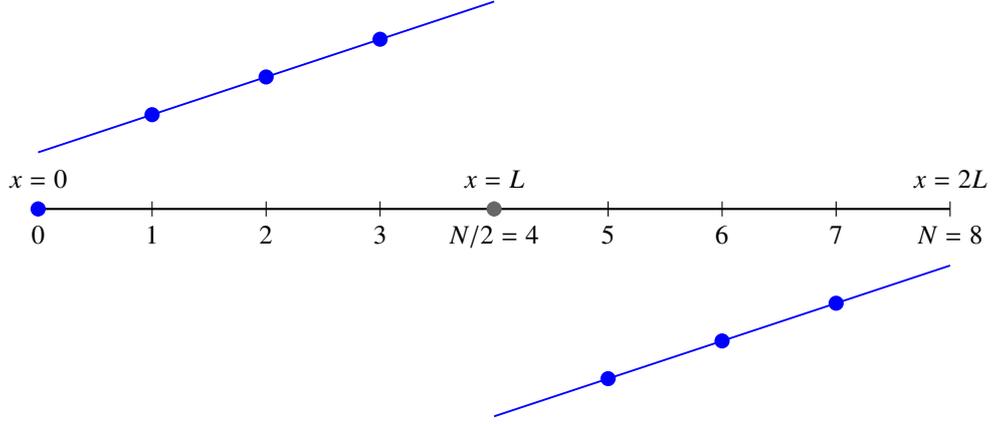

%
We consider on the domain~$\Omega = (0, \, L) \subset \mathbb R^1$ the Poisson equation 
\begin{equation} \label{eq:poisson_first}
	- \frac{\D^2 u(x)}{\D x^2} = f(x), \quad \forall x \in \Omega \, ,
\end{equation}
with the homogeneous Dirichlet boundary conditions
\begin{equation}
	u(0) = u(L) = 0 \, .
\end{equation}
We extend the source and solution fields,~$f(x)$ and~$u(x)$, to the domain $ \Omega_{\text E} = (0, \, 2L) $ as follows
\begin{subequations}
\begin{align} 
	f_{\text E}(x) =  
	\begin{cases}
		0  & \text{if} \quad x = 0, \, L , \, 2L, \\
		f(x) & \text{if} \quad x \in (0, \, L) \, , \\
		-f(2L-x) & \text{if} \quad x \in (L, \,  2L) \, .
	\end{cases} \label{eq:source_reflection}
\\
	u_{\text E}(x) = 
	\begin{cases}
		0  & \text{if} \quad x = 0, \, L , \, 2L, \\
		u(x) & \text{if} \quad x \in (0, \, L) \, , \\
		-u(2L-x) & \text{if} \quad x \in (L, \,  2L) \, ;
	\end{cases}
\end{align}
\end{subequations}
see Figure~\ref{fig:extended_domain}. It will be seen that this antisymmetric extension ensures that $ u(x) $ satisfies the Dirichlet boundary conditions, making it compatible with a Fourier sine series expansion. Note that the symmetry properties of the solution $ u(x) $ depend on the order of the differential equation: for even-order differential equations, the source term $ f(x) $ and the solution $ u(x) $ share the same symmetry; for odd-order differential equations, the symmetries of $ u(x) $ and $ f(x) $ are opposite (i.e., one is symmetric while the other is antisymmetric).

On the extended domain~$\Omega_{\text E}$, the solution and source are approximated using truncated complex Fourier series
\begin{equation} \label{eq:exp_series_uf}	
	u_{\text E}(x) \approx u_\text{E}^h(x)= \frac{1}{\sqrt{N}} \sum_{k=-N/2}^{N/2-1} \hat u_{ \text E,  k}  e^{ i \frac{ \pi k x}{L}} \, , \qquad f_{\text E}(x) \approx f_{\text E}^h(x) = \frac{1}{\sqrt{N}} \sum_{k=-N/2}^{N/2-1} \hat f_{\text E, k} e^{ i \frac{ \pi k x}{L}}  \, ,
\end{equation}
where $  \hat u_{ \text E,  k} $ and $ \hat f_{\text E, k}  $ are Fourier coefficients, $ i^2 = -1 $, and $e^{ i \pi k x/ L} = \cos( \pi k x / L) +i \sin( \pi k x / L)$. To render the index~$k$ non-negative,  we  introduce the relabelling 
\begin{equation} 
	\label{eq:relabelling}
	r(k) = \begin{cases}
		k & 0 \le k < N/2 \\
		k - N & N/2 \le k < N \, ;
	\end{cases}
\end{equation}
so that
\begin{equation} 
	\label{eq:components_to_frequencies}
	u_{\text E}^h(x) =  \frac{1}{\sqrt N} \sum_{k=0}^{N-1} \hat u_{\text E, k} e^{ i \xi_k x }  \, , \qquad f_{\text E}^h(x) = \frac{1}{\sqrt N} \sum_{k=0}^{N-1} \hat f_{\text E, k} e^{ i \xi_k x }  \, , \qquad \text{where} \; \,  \xi_k = \frac{ \pi r(k)}{L} \, .
\end{equation}

Substituting the Fourier basis functions into the governing equation~\eqref{eq:poisson_first} yields 
\begin{equation}
	\sum_{k=0}^{N-1} \left (  \xi_k^2 \hat u_{\text E, k}  - \hat f_{\text E, k}  \right ) e^{i \xi_k x}  = 0 \, , 
\end{equation}
which is required to be zero in a weak sense with the Fourier basis functions as the test functions, i.e., 
\begin{equation}
	\mathlarger{\int}_{\Omega_{\text E}} \left [  \sum_{k=0}^{N-1} \left (  \xi_k^2 \hat u_{\text E, k}  - \hat f_{\text E, k}  \right )  e^{i \xi_k x} \right ]  e^{-i \xi_j x} \D x = 0 \, .
\end{equation}
Given the orthogonality of the Fourier basis, this yields the algebraic relation
\begin{equation} \label{eq:discrete_poisson}
	\hat u_{\text E, k} =   \frac{\hat f_{\text E, k}}{\xi_k^2} \, , \qquad \forall  k \neq 0 \, .
\end{equation}
The wavenumber~$k=0$ can be omitted because~$f_{\text E} (x)$ has zero mean, and therefore~$\hat f_{\text E, 0}=0$. Introducing the Fourier domain solution in the first equation in~\eqref{eq:components_to_frequencies} yields 
\begin{equation} \label{eq:discrete_poisson_rs}
	u_{\text E}^h(x) =  \frac{1}{\sqrt N} \sum_{k=0}^{N-1} \frac{\hat f_{\text E, k}}{\xi^2_k}  e^{ i \xi_k x }.
\end{equation}
This solution is inherently periodic and exhibits antisymmetry about the centre of the domain.
%
\subsection{Solution by discrete Fourier transform}
\label{sec:dft-soln}
%
We determine the Fourier coefficients $ \hat{f}_{\text E, k}$ corresponding to the source field approximation $ f^h_{\text E}(x) $ in~\eqref{eq:exp_series_uf} using the discrete Fourier transform (DFT). To achieve this, we discretise the extended domain $ \Omega_\text{E} $ with $N$ cells and equidistant grid points, as illustrated in Figure~\ref{fig:extended_domain}. Thus, the physical domain $ \Omega$ contains $ N/2 $ grid points and $ N/2 $ grid cells.  The grid point coordinates are defined as  
\begin{equation}
	x_k = k \frac{2L}{N}, \quad k = 0, 1, \dotsc, N-1 \, .
\end{equation}  
The source field values at the grid points are collected into the antisymmetrically extended vector  
\begin{equation} \label{eq:source_vector}
	\vec{f}_{\text E} = \begin{pmatrix} 0 &  f(x_1) &  \dotsc &  f(x_{N/2-1}) &  0 &  -f(x_{N/2-1}) &   \dotsc &   -f(x_1) \end{pmatrix}^\trans \, .
\end{equation}  
The components corresponding to the first and centre grid points are zero as required by~\eqref{eq:source_reflection}. We now introduce the antisymmetric extension matrix~$\vec R \in \mathbb R^{N \times N/2}$ to determine the reflected source vector~\mbox{$\vec f_\text E \in \mathbb R^{N}$} from the source vector~\mbox{$\vec f \in \mathbb R^{N/2}$} containing grid point values only associated to the physical domain, i.e.,
\begin{equation} \label{eq:reflection_source_vector}
	\vec f_{\text E} = \vec R \vec f \, , 
\end{equation}
where
\begin{equation} \label{eq:source_vector_half}
	\vec f  = \begin{pmatrix} 0 &  f(x_1) &  \dotsc & f(x_{N/2-1})  \end{pmatrix}^\trans  \, ,  
\end{equation}
and
\begin{equation} \label{eq:reflection_matrix}
	\vec R =
	\begin{pmatrix*}[r]
		0 & 0 & 0 & \cdots & 0 \\
		0 & 1 & 0 &  \cdots & 0  \\
		0 & 0 & 1 &  \cdots & 0 \\
		\vdots & \vdots &\vdots &  \ddots &  \vdots  \\
		0 & 0 &  0  & \cdots & 1  \\
		0 & 0 &  0  & \cdots & 0  \\
  		0 & 0 &  0  & \cdots & -1  \\
		\vdots & \vdots &\vdots &  \iddots & \vdots   \\
		0 & 0 &  -1 & \cdots & 0 \\
		0 & -1  & 0 & \cdots & 0 \\
	\end{pmatrix*} \, .
\end{equation}
The so-defined reflection matrix~$\vec R$ ensures that the source vector~$\vec f_{\text E}$ is antisymmetric.

The DFT of the source vector~$\vec f_{\text E}$ is defined as  
\begin{equation} \label{eq:dft}
	\hat{\vec f}_\text{E} = {\vec F_N}{\vec f}_{\text E} \, ,
\end{equation}
with the Fourier matrix 
\begin{equation}
	{\vec F}_N = 
	\frac{1}{\sqrt {N}}  {\left (  \omega_N^{jk}  \right )}_{j,k=0, \dotsc, N-1 } \, ,
\end{equation}
and the $N$-th root of unity $\omega_N = e^{i 2 \pi /N}$.
The Fourier matrix is unitary, i.e. $ \vec F_N \vec F^\dagger_N = \vec I$, where the dagger $\dagger$ denotes the conjugate transpose, and can be written as 
\begin{align} \label{eq:ft_matrix_explicit}
	{\vec F}_N = \frac{1}{\sqrt N}  
	\begin{pmatrix}
		1 & 1 & 1 & 1 &  \cdots & 1 &  \cdots & 1 & 1 & 1 \\
		1 & \omega_N & \omega_N^2 & \omega_N^3 & \cdots &  \omega_N^{N/2} & \cdots &  \omega_N^{-3} & \omega_N^{-2} & \omega_N^{-1} \\[0.15em]
		1 & \omega_N^2 & \omega_N^4 & \omega_N^{-2}& \cdots  & 1 & \cdots  & \omega_N^{2} &  \omega_N^{4} & \omega_N^{-2}  \\[0.15em]
		\vdots & \vdots & \vdots & \vdots & \ddots & \vdots  &  \iddots & \vdots & \vdots & \vdots  \\[0.1em]
		1 & \omega_N^{N/2} & 1 & \omega_N^{N/2} & \cdots & \omega_N^{N/2} & \cdots & \omega_N^{N/2} & 1 &   \omega_N^{N/2}  \\[0.15em]
		\vdots & \vdots & \vdots & \vdots &  \iddots & \vdots  & \ddots & \vdots & \vdots & \vdots  \\[0.1em]
		1 & \omega_N^{-2} & \omega_N^{4} & \omega_N^{2} & \cdots & 1 &\cdots  & \omega_N^{-2} & \omega_N^{-4} & \omega_N^{2}  \\[0.15em]
		1 & \omega_N^{-1} & \omega_N^{-2} & \omega_N^{-3} & \cdots & \omega_N^{N/2} & \cdots &  \omega_N^{3} & \omega_N^{2}  & \omega_N 
	\end{pmatrix} \, , 
\end{align} 
taking into account~$\omega_N^k = \omega_N^{(k \mod N)}$ and  $\omega_N^{k} = \omega_N^{-(N-k)}$. It is easy to verify that the product~$\vec F_N \vec R$ yields a matrix with only sine terms, with each column representing one of the Fourier basis functions.  The cosine terms in the original~$\vec F_N$ are eliminated by multiplication with the reflection matrix~$\vec R$. 
The first~$N/2$ rows correspond to grid points within the domain of interest~$\Omega$. 

With the help of definitions of~$\vec f_\text{E}$ and DFT, \eqref{eq:reflection_source_vector} and~\eqref{eq:dft},  the solution~\eqref{eq:discrete_poisson_rs} in the Fourier domain can be compactly written as
\begin{equation} \label{eq:solution_extended}
	\vec u_{\text E} =     \vec F_N^\dagger \vec D_{\text E}^{-1} \vec F_N \vec R \vec f \, , 
\end{equation}
where~$\vec D_{\text E} \in \mathbb R^{N\times N}$ is a diagonal matrix with the entries
\begin{equation}
	\vec D_{\text E} = \left ( \frac{\pi}{L} \right )^2  \diag \begin{pmatrix} 1 & 1 & 2^2 & 3^2 & \cdots & (N/2)^2 &  \cdots & 3^2  &  2^2 & 1 \end{pmatrix} \, .
\end{equation}
We can rewrite the solution equation such that 
\begin{equation}
	\vec u =   \vec R^\trans \vec F_N^\dagger \vec R  \vec D^{-1} \vec R^\trans \vec F_N \vec R \vec f \, , 
\end{equation} 
where 
\begin{equation} \label{eq:solution_poisson_fsp}
	\vec D  = \left ( \frac{\pi}{L} \right )^2  \diag \begin{pmatrix} 1 & 1 & 2^2 & 3^2 & \cdots & (N/2-1)^2 \end{pmatrix} \, .
\end{equation}
In comparison to~\eqref{eq:solution_extended}, this has the advantage that the size of the matrix~$\vec D \in \mathbb R^{N/2 \times N/2}$ that must be quantum encoded later is halved. Furthermore, we note that the expression~$ \vec R^\trans \vec F_N \vec R$ is known as the discrete sine transform (DST) of type~I~\cite{strang1999discrete}. 

%
\subsection{Other boundary value problems\label{sec:other-PDE}}
%
The introduced Fourier spectral method can easily be generalised. Because we use trigonometric sine functions as basis functions, the occurring derivatives in the differential equation must be all of even order (i.e., zeroth, second, fourth, ...).   We consider as an example problem on~$\Omega  = (0, \, L) \subset \mathbb R^1$ the fractional differential equation  
\begin{equation}  \label{eq:fractional_ode}
	\left ( \kappa^2  - \frac{\D^2 }{\D x^2} \right )^\beta  u ( x) = \frac{1}{\tau} f (x) \, ,  \qquad \forall  x \in \Omega \, ,
\end{equation}
with the homogeneous Dirichlet boundary conditions 
\begin{equation}
	u(0) = u(L) = 0 \, , 
\end{equation}
and the three prescribed parameters~$\kappa, \beta, \tau \in \mathbb R$. 

As in the Poisson-Dirichlet problem, first the source and solution fields~$f(x)$ and~$u(x)$  are antisymmetrically extended to the larger domain~$\Omega_\text{E} = (0, 2L)$. Subsequently, the obtained fields~$f_{\text E}(x)$ and~$u_{\text E}(x)$ are approximated using a truncated Fourier series. A straightforward derivation shows that the Fourier domain solution satisfies the algebraic equation
\begin{equation}
	\hat u_{\text E, k} = \frac{\hat f_{\text E, k} }{\tau (\kappa^2 +   \xi_k^2  )^\beta  }  \, ,  \qquad \forall k \neq 0 \, .
\end{equation}
Hence, the solution of the fractional differential equation in the extended domain is given by 
\begin{equation}
	u_{\text E}^h(x) =  \frac{1}{\sqrt N} \sum_{k=0}^{N-1}  \frac{\hat f_{\text E, k} }{\tau (\kappa^2 +   \xi_k^2  )^\beta  }  e^{i \xi_k x} \, .
\end{equation}
Again,~$\hat u_{\text E, k} $ and~$\hat f_{\text E, k}$ coincide with the discrete Fourier transforms of~$u _{\text E}(x)$ and~$f_{\text E}(x)$ sampled over an equidistant grid.  By retracing the same steps as in the Poisson-Dirichlet case, we obtain the discrete solution vector
\begin{equation} \label{eq:solution_fractional}
	\vec u  = \vec R^\trans \vec F^\dagger \vec R {\vec D}^{-1} \vec R^\trans \vec F \vec R \vec f \, , 
\end{equation} 
where 
\begin{equation} \label{eq:solution_poisson_fractional}
	\vec D  = \tau  \diag \begin{pmatrix} \kappa^2 & \kappa^2+ \left ( \frac{\pi}{L} \right )^2 &  \kappa^2+ \left ( \frac{2 \pi}{L} \right )^2  &   \kappa^2+ \left ( \frac{3 \pi}{L} \right ) ^2  & \cdots &  \kappa^2+ \left ( \frac{(N/2-1) \pi}{L} \right ) ^2   \end{pmatrix}^\beta \, .
\end{equation}
%

%
\subsection{Arbitrary Dirichlet boundary conditions}
%
The considered boundary value problems are all linear, which allows us to use the superposition principle when applying inhomogeneous Dirichlet boundary conditions.  The main idea is to choose first a sufficiently smooth function~$g(x)$ that satisfies the prescribed Dirichlet boundary conditions, but otherwise has arbitrary values within the domain. For instance, in the case of the Poisson-Dirichlet problem, we have 
\begin{equation}
\begin{aligned} \label{eq:poisson_boundary}
	& - \frac{\D^2 u(x)}{\D x^2} = f(x), \quad  \forall x \in (0, \, L) \, , \\
	& u(0) = g(0) \, , \quad  u(L) = g(L)  \, . 
\end{aligned}
\end{equation}
The function $ g (x)$ can, for instance, be constructed by interpolating the prescribed Dirichlet boundary conditions using Lagrange basis functions. We can now define a solution field~$v(x) = u(x) - g(x)$ that satisfies 
\begin{equation}
\begin{aligned} 
	& - \frac{\D^2 v(x)}{\D x^2} = f(x) +  \frac{\D^2 g(x)}{\D x^2}  , \quad  \forall x \in (0, \, L) \, , \\
	& v(0) = 0 \, , \quad  v(L) = 0  \, . 
\end{aligned}
\end{equation}
Consequently, in the Fourier spectral method, arbitrary Dirichlet boundary conditions can be accommodated by suitably modifying the source field of a problem with homogeneous Dirichlet boundary conditions. 
%

%

%
%
\section{Quantum Fourier spectral method}
%
We introduce next the quantum implementation of the Fourier spectral method for arbitrary Dirichlet boundary conditions. The forcing vector is encoded as the initial state of a quantum system, and the partial differential equation is solved by applying a sequence of unitary transformations to the state. In this section, we switch to Dirac notation, which is analogous to the standard vector notation, to distinguish between quantum mechanical and classical vectors. Introductions to Dirac notation can be found in standard texts~\cite{ikeAndMike,kaye2006introduction,rieffel2011quantum} and our earlier work~\cite{liu2024towards}.  
%
\subsection{Quantum Fourier transform (QFT)}
%
In Dirac notation, the Fourier transform~\eqref{eq:dft} is expressed as 
\begin{equation}
	\widehat{\ket f}  = F_N \ket f \, ,
\end{equation}
where~$\widehat{\ket f}, \, \ket f \in \mathbb C^N$ are the $N=2^n$ dimensional state vectors of a quantum system with~$n$ qubits, before and after the application of the Fourier matrix~$F_N \in \mathbb C^{N \times N}$. The two vectors in componentwise form are given by
\begin{equation}
	\ket f = \begin{pmatrix}f_0 \\ f_1 \\ \vdots \\ f_{N-1} \end{pmatrix} =   \sum_{k=0}^{N-1}  f_k   \ket k  \, , \quad \widehat{\ket f} =  \begin{pmatrix} \hat f_0 \\ \hat f_1 \\ \vdots \\ \hat f_{N-1} \end{pmatrix} = \sum_{k=0}^{N-1} \hat f_k   \ket k \, , 
\end{equation}
where~$\ket k$, with $k \in \{0, \, 1,  \dotsc, \, N-1 \}$,  are the standard basis vectors, i.e., 
\begin{equation}
	\ket 0  =  \begin{pmatrix}1 \\ 0 \\ \vdots \\ 0 \end{pmatrix} \, ,  \quad  \ket 1 = \begin{pmatrix}0 \\ 1 \\ \vdots \\ 0 \end{pmatrix} \, , \quad \dotsc \, , \quad   \ket {N-1}  =  \begin{pmatrix}0 \\ 0 \\ \vdots \\ 1 \end{pmatrix} \, ,
\end{equation}
and~$f_k$ and~$ \hat f_k$  are the respective coefficients. In quantum computing, the standard basis is often called the computational basis, and the coefficients are called amplitudes. The basis vectors can be indexed using binary numbers by noting that 
\begin{equation}
k = k_0 2^{n-1} + \dotsc + k_{n-2} 2^1 + k_{n-1} 2^0  = \sum_{j=0}^{n-1} k_j 2^{n-1-j}\, \, \text{ and } \, \,  k_0, \,   \dotsc , \, k_{n-2}, \, k_{n-1}  \in \{ 0, \, 1\} \, ,
\end{equation}
 and writing~$\ket k \equiv \ket{k_0 \cdots  k_{n-2} k_{n-1}} \equiv \ket{k_0} \cdots \ket{k_{n-2}} \ket{k_{n-1}}$.
  
\begin{figure}[]
	\centering
	\scalebox{.88}{
		\input{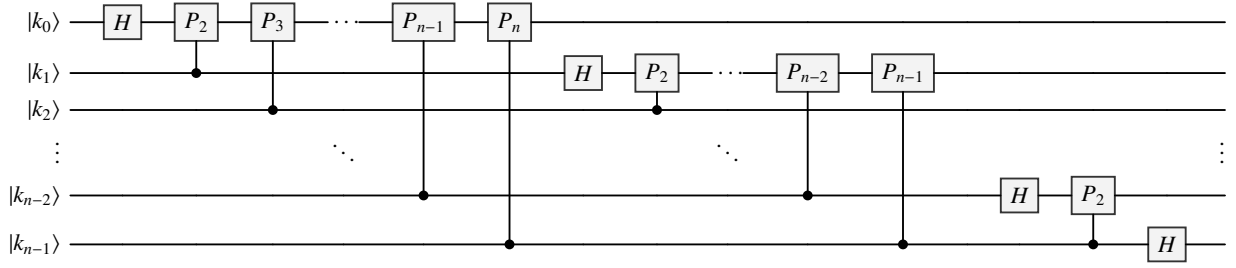}
	}
	\caption{Quantum circuit for~QFT. An input vector~$\sum_{k=0}^{N-1}f_k \ket k$ is mapped to the output vector~$\sum_{k=0}^{N-1} \hat f_k \ket k $ (pursuant to some relabelling of the components of the output vector). \label{fig:qft-circ}}
\end{figure}

The Fourier matrix is inherently unitary,~i.e. $F_N F_N^\dagger=I$, as its row and column vectors are mutually orthonormal. $F_N$ is implemented by multiplicatively decomposing it into~$2\times2$ and~$4\times4$ unitary matrices, referred to as gates.  It bears emphasis here that any product of unitary matrices is again a unitary matrix. In coming up with a suitable decomposition, it is convenient to consider the action~$F_N$ on a single basis vector $\ket {k}$,  yielding the  column~$k$ of~$F_N$.   Subsequently, the derived decomposition can be applied to the entire vector~$\ket f$ by making use of quantum superposition. The key realisation in implementing~$F_N$ is that its column~$k$ can be expressed as the tensor product, i.e. Kronecker product, of~$n$ two-dimensional vectors, 
\begin{equation} \label{eq:qft_product}
\begin{aligned} 
	 \frac{1}{\sqrt N} \otimes_{j=0}^{n-1} ( \ket 0 +  \omega_N^{ 2^{n-1-j}  k }\ket 1 ) &= \frac{1}{\sqrt N}   \begin{pmatrix}1  \\ \omega_N^{ 2^{n-1}  k } \end{pmatrix} \otimes  \cdots  \otimes \begin{pmatrix}1  \\ \omega_N^{ 2^{2}  k } \end{pmatrix}  \otimes \begin{pmatrix}1  \\ \omega_N^{ 2^{1}  k } \end{pmatrix}  \otimes \begin{pmatrix}1  \\ \omega_N^{ 2^{0}  k } \end{pmatrix}  \\
	 &= \frac{1}{\sqrt N}  	 \begin{pmatrix} 1 & \omega_N^k & \omega_N^{2 k} & \omega_N^{3k} & \omega_N^{4k} & \dotsc &  \omega_N^{(2^{n}-1)k }	 \end{pmatrix}^\trans  \, ,\\ 
\end{aligned}
\end{equation}
cf.~\eqref{eq:ft_matrix_explicit}.  Here, each two-dimensional vector represents the state vector of a single qubit. The first line in~\eqref{eq:qft_product} is rewritten by expressing the index~$k$ as a binary, i.e., 
\begin{equation} \label{eq:tensor_product_qft}
 \frac{1}{\sqrt N}   \begin{pmatrix}1  \\ e^{ 2 \pi i  \left ( \frac{k_{n-1}}{2} \right )} \end{pmatrix} \otimes  \cdots  \otimes \begin{pmatrix}1  \\ e^{ 2 \pi  i  \left ( \frac{k_2}{2}   + \dotsc + \frac{k_{n-1}}{2^{n-2}}  \right )} \end{pmatrix}  \otimes \begin{pmatrix}1  \\ e^{ 2 \pi i   \left ( \frac{k_1}{2}   + \dotsc + \frac{k_{n-1}}{2^{n-1}}  \right )} \end{pmatrix}  \otimes \begin{pmatrix}1  \\ e^{2 \pi i  \left ( \frac{k_0}{2}   + \dotsc + \frac{k_{n-1}}{2^n}  \right )} \end{pmatrix}  \, .
\end{equation}
The identity~$e^{2\pi i} \equiv 1$ and its periodicity are used in deriving this expression. Hence, the state vector of the entire system consists of the tensor product of the states of~$n$ qubits and depends on the binary representation of the column index~$k$ of~$F_N$.  As depicted in the circuit diagram in Figure~\ref{fig:qft-circ}, the state vector~\eqref{eq:tensor_product_qft} can be constructed by using only Hadamard and controlled phase gates defined as
\begin{equation}
	H = \frac{1}{\sqrt{N}} \begin{pmatrix*}[r] 1 & 1 \\ 1 & -1 \end{pmatrix*} \, , \quad P_l  = \begin{pmatrix} 1 & 0 \\ 0  & e^{2\pi i / 2^l}  \end{pmatrix}  \, .
\end{equation}
As an example, consider the case~$k=0$ with the input state vector~$\ket{0 \dotsc 000}$ yielding~$F_N$'s first column  
\begin{equation} 
 \frac{1}{\sqrt N}   \begin{pmatrix}1  \\ 1  \end{pmatrix} \otimes  \cdots  \otimes \begin{pmatrix}1  \\ 1 \end{pmatrix}  \otimes \begin{pmatrix}1  \\  1\end{pmatrix}  \otimes \begin{pmatrix}1  \\1  \end{pmatrix}  =
 \frac{1}{\sqrt N} 
 \begin{pmatrix}
  1 & 1 & 1 & 1 &1 & \dotsc &  1   
 \end{pmatrix}^\trans  \, .
\end{equation}
This state is constructed by applying Hadamard gates to all qubits. None of the controlled phase gates are applied because the respective control qubits are in state~$\ket 0$. On the other hand, when~$k=1$ with the input state vector~$\ket{0 \dotsc 001}$ the state vector that must be constructed reads  
\begin{equation} \label{eq:tensor_product_qft_k1}
 \frac{1}{\sqrt N}   \begin{pmatrix}1  \\  1 \end{pmatrix} \otimes  \cdots  \otimes \begin{pmatrix}1  \\ e^{\frac{2 \pi i }{2^{n-2}} } \end{pmatrix}  \otimes \begin{pmatrix}1  \\ e^{\frac{2 \pi i}{2^{n-1}}}      \end{pmatrix}   \otimes \begin{pmatrix}1  \\ e^{\frac{ 2 \pi i}{2^{n}}}      \end{pmatrix} 
 = \frac{1}{\sqrt N} 
 \begin{pmatrix}
 	1  & e^\frac{2 \pi i}{2^n}  & e^\frac{4 \pi i}{2^n}  & e^\frac{6 \pi i}{2^n} & e^\frac{8 \pi i}{2^n} & \dotsc  & e^\frac{2(2^n-1) \pi i}{2^n}
 \end{pmatrix}^\trans \, .
\end{equation}
This state is constructed by applying controlled phase gates~$P_n, \,  P_{n-1} ,\, \dotsc , \, P_2 $ conditioned on~$\ket{k_{n-1}}$ is in state~$\ket 1$.  In the circuit in Figure~\ref{fig:qft-circ}, the state~\eqref{eq:tensor_product_qft_k1} is first generated in reverse order and must be subsequently corrected using swap gates mapping~$\ket{ k_{n-1} k_{n-2} k_{n-3}  \cdots k_0 }$ to~$\ket{k_0 \cdots  k_{n-3} k_{n-2} k_{n-1}}$. These swap gates are omitted in Figure~\ref{fig:qft-circ}. The quantum circuit takes  as an input the state vector~$\sum_{k=0}^{N-1} f_k \ket k$ and returns as an output the state vector~$\sum_{k=0}^{N-1} \hat f_k \ket k$. Hence, the entire quantum state is transformed in one fell swoop owing to quantum parallelism.

In passing, we note that the QFT algorithm has the complexity $O(n^2)$ where~$n = \log N$~\cite[Ch. 5]{ikeAndMike}. In comparison, the complexity of fast Fourier transform (FFT) expressed in terms of~$n$ is~$O(n 2^n)$. That is, QFT is exponentially faster than FFT.

%
\subsection{Unitary reflection matrix \label{sec:quantum-reflect}}
%
We consider next the quantum implementation of the non-unitary reflection matrix~$\vec R  \in\mathbb R^{N/2\times N}$ introduced in~\eqref{eq:reflection_matrix}. To this end, the matrix~$\vec R$ is embedded, or in quantum computing terminology block-encoded, in the unitary matrix 
\begin{equation}
	U_R = \frac{1}{\sqrt 2 }
	\begin{pmatrix}
		0 & 0 & 0 & \cdots & 0 & \sqrt 2 &0 &0 & \cdots & 0  \\ 
		0 & i & 0 & \cdots & 0 & 0 & 1 & 0 & \cdots & 0  \\
		0 & 0 & i & \cdots & 0 & 0 & 0 & 1 & \cdots & 0  \\
		\vdots & \vdots & \vdots & \ddots & \vdots & \vdots& \vdots & \vdots &  \ddots &\vdots \\
		 0 & 0 & 0 & \cdots & i & 0& 0& 0& \cdots & 1 \\
		 \sqrt 2 & 0 & 0 & \cdots & 0  & 0 & 0 & 0 & \cdots & 0  \\
		 0& 0& 0& \cdots & -i  & 0& 0 & 0 & \cdots & 1  \\
		\vdots & \vdots & \vdots & \iddots & \vdots & \vdots & \vdots & \vdots & \iddots & \vdots \\
		 0 & 0 & -i & \cdots & 0 & 0 & 0 & 1 & \cdots & 0 \\
		 0 & -i & 0 & \cdots & 0 & 0 & 1 & 0 & \cdots & 0 
	\end{pmatrix} \, .
	\label{eq:matrixT}
\end{equation}
Except for the first column, its left half is equal to~$i \vec R$ up to a scaling factor. The first column and the second half of~$U_R$ are chosen so that all rows and columns are mutually orthonormal, as necessary for unitarity. Although the first column is different from~$\vec R$'s first column, it is inconsequential for the final solution of a zero Dirichlet boundary value problem.  

\begin{figure}[]
	\centering
	\input{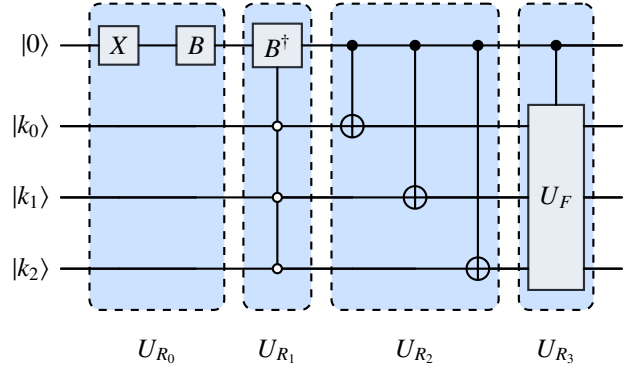}
	\caption[]{Quantum circuit for the reflection unitary $U_R= U_{R_3} U_{R_2} U_{R_1} U_{R_0}$. An input vector~$\sum_{k=0}^{7} f_k \ket 0 \ket k$ is mapped to the antisymmetrically extended output vector~$\sum_{k=0}^{15} f_k  \ket k$.  \label{fig:matT-mcx}}
\end{figure}

The unitary $U_R$ is again implemented by decomposing it into elementary gates, i.e.~$2\times2$ and~$4\times4$ unitary matrices. A decomposition of~$U_R$ following~\cite{klappenecker2001discrete} is shown in Figure~\ref{fig:matT-mcx}; alternative decompositions are possible. For concreteness, we choose for the extended domain a discretisation with~$N=16$ grid points, i.e. four qubits. Consequently, the problem domain is discretised with~$8$ grid points and the respective field vector values are encoded as the state vector of the three qubits~$\ket{k_0 k_1 k_2}$. 
The depicted decomposition can be understood as the product of the four unitaries~$U_{R_0}$, $U_{R_1}$, $U_{R_2}$, $U_{R_3}$, each of which will be discussed in the following.  

The leftmost unitary~$U_{R_0}$ in the circuit diagram is given by
\begin{equation} \label{eq:U_R_0}
	U_{R_0} = (B X) \otimes I^{\otimes 3}  =  ( H S X) \otimes I^{\otimes 3} =  \frac{1}{\sqrt 2} \left ( \begin{pmatrix*}[r] 1 & 1 \\ 1 & -1\end{pmatrix*}  \begin{pmatrix} 1 & 0 \\ 0 & i\end{pmatrix}   \begin{pmatrix} 0 & 1 \\ 1 & 0\end{pmatrix}   \right )  \otimes I^{\otimes 3} =  \frac{1}{\sqrt 2} \begin{pmatrix*}[r] i I ^{\otimes 3} &  I ^{\otimes 3}    \\  -i I^{\otimes 3}   & I^{ \otimes 3} \end{pmatrix*} \, ,
\end{equation}
where~$I  \in \mathbb R^{2 \times 2}$ and $I^{\otimes 3} \in \mathbb R^{8 \times 8}$  are identity matrices, and~$X$, $H$ and~$S$ are the single-qubit Pauli~$X$, Hadamard and phase gates which are available in most quantum computing platforms.   The unitary~$U_{R_0} \in \mathbb C^{16 \times 16}$ consists of the four blocks indicated in~\eqref{eq:U_R_0}. Notice that the identity gates in~\eqref{eq:U_R_0} are, as usual, omitted in the circuit diagram. 

The unitary~$U_{R_1}$ consists of the conjugate transpose~$B^\dagger$ conditioned on the state~$\ket{k_0k_1k_2}$. Specifically, $B^\dagger$ is only applied to the ancilla qubit when~$\ket{k_0k_1k_2} \equiv \ket {000} $ as indicated by the empty circles. Consequently, it applies a unitary transformation in the subspace spanned by the two basis vectors~$\ket 0 \ket{000}$ and~$\ket 1 \ket{000}$,  else the identity gate~$I$ is applied. Based on the circuit diagram, we can write 
\begin{equation}
\begin{split}
	U_{R_1} &= B^\dagger \otimes \left (  \ket 0 \bra 0  \right )^{\otimes 3} +  \left ( I \otimes  \displaystyle \sum_{\substack{k_0=0, k_1=0, k_2 =0 \\ k_0 = k_1 = k_2 \neq 0 }}^{1}  \ket {k_0} \bra {k_0} \otimes \ket {k_1} \bra {k_1} \otimes \ket {k_2} \bra {k_2}  \right )  \, , 
 \end{split}
\end{equation}
and, more explicitly,  
\begin{equation}
	U_{R_1} =  \begin{pmatrix}	\frac{1}{\sqrt 2}   & 0 & \cdots & 0 &    	\frac{1}{\sqrt 2}   & 0 & \cdots & 0 \\ 
							  0 & 1 & \cdots & 0 & 0 & 0 & \cdots & 0  \\ 
							  \vdots & \vdots & \cdots & \vdots & \vdots & \vdots & \cdots & \vdots  \\ 
							  0 & 0 & \cdots & 1 & 0 & 0 & \cdots & 0  \\ 
							  \frac{- i}{\sqrt 2}  & 0 & \cdots & 0 & \frac{i}{\sqrt 2} & 0 & \cdots & 0  \\
							  0 & 0 & \cdots & 0 & 0 & 1 & \cdots & 0  \\ 
							  \vdots & \vdots & \cdots & \vdots & \vdots & \vdots & \cdots & \vdots  \\ 
							  0 & 0 & \cdots & 0 & 0 & 0 & \cdots & 1  \\ 
							   \end{pmatrix} \, .
\end{equation}

The next unitary~$U_{R_2}$ is composed of three two-qubit $CNOT$ gates, which are applied only when the top (ancilla) qubit is in state~$\ket 1$, as indicated by the filled circles. If the ancilla qubit is in state~$\ket 0$ the identity gate~$I$ is applied.  A single~$CNOT$ gate is a controlled~$X$ gate.  Hence,~$U_{R_2}$ takes the form 
\begin{align}
	U_{R_2} & =  \ \ket 0 \bra 0 \otimes I^{\otimes 3} +  \ket 1 \bra 1  \otimes X^{\otimes 3}  =  \begin{pmatrix} 1 & 0 \\ 0 & 0 \end{pmatrix}  \otimes \begin{pmatrix} 1 & 0 \\ 0 & 1 \end{pmatrix}^{\otimes 3} +  \begin{pmatrix} 0 & 0 \\ 0 &  1\end{pmatrix}  \otimes \begin{pmatrix} 0 & 1 \\ 1 & 0 \end{pmatrix}^{\otimes 3} = \begin{pmatrix} I^{\otimes 3} & 0 \\ 0 & X^{\otimes 3}\end{pmatrix} \, .
\end{align}
The multiplication of the unitaries~$ U_{R_1}U_{R_0}$ defined so far with~$U_{R_2}$ from the left has the effect of moving the diagonal entries in the two blocks on the lower half of~$ U_{R_1}U_{R_0}$ along the anti-diagonal. It bears emphasis that the order of the unitaries in the circuit diagram and matrix multiplication is reversed. The circuit diagram is read from left to right, whereas the matrix product is defined right to left. 

\begin{figure}[]
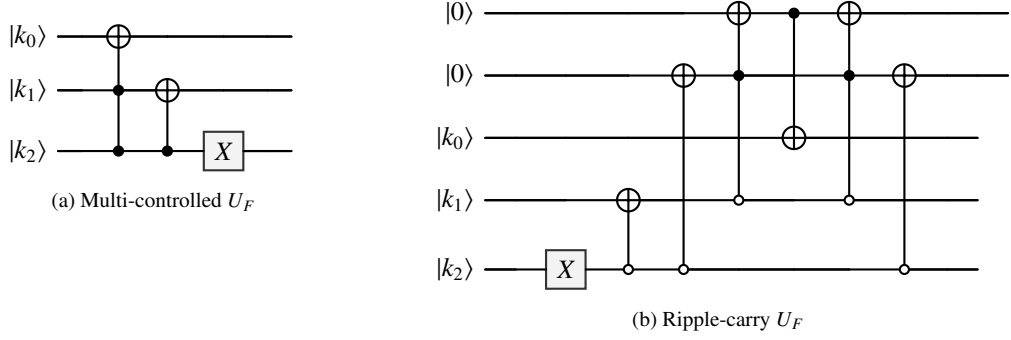

	\centering
	\subfloat[][Multi-controlled $U_F$\label{fig:matp_mcx}]{
		\input{figures/circuits/matP_mcx_q3}
	}
	\hspace{0.075\textwidth}
	\subfloat[][Ripple-carry $U_F$\label{fig:matp_cascade}]{
		\raisebox{-0.3\height}{
			\input{figures/circuits/matP_cascade_q3}}
	}
	\caption[]{Two alternative quantum circuits for the forward shift unitary~$U_F$,  (a) using a cascade of controlled~$X$ gates and (b) based on ripple-carry addition. In (b), the top two qubits are ancilla qubits. Both circuits map an input vector~$\sum_{k=0}^{7} f_k \ket k $ to an output vector $\sum_{k=0}^{7} \ket {(k+1) \mod 2^3}$. The ancilla qubits in (b) have the value~$\ket 0$ before and after execution of the circuit.\label{fig:matp}}
\end{figure}

The final unitary~$U_{R_3}$ in the circuit diagram is the controlled (periodic) forward shift unitary, which applies a down shift by one to the entries in the lower half of~$U_{R_2} U_{R_1}U_{R_0}$. It follows from the circuit diagram that
\begin{equation}
	U_{R_3} = \ket 0 \bra 0 \otimes I^{\otimes 3} + \ket 1 \bra 1 \otimes U_F \, ,
\end{equation}
where the forward shift unitary is defined as~$U_F \colon \ket k \mapsto \ket {(k+1) \mod 2^3}$, i.e., 
\begin{equation}
U_F = \begin{pmatrix} 0 & 0 & \cdots & 0 & 1 \\
1 & 0 & \cdots & 0 & 0 \\
0 & 1 & \cdots & 0 & 0 \\
\vdots & \vdots & \cdots & \vdots & \vdots \\
0 & 0 & \cdots & 1 & 0 
\end{pmatrix} \, .
\end{equation}
There are several possible ways to implement the unitary~$U_F$. In Figure~\ref{fig:matp_mcx}, we present an implementation of~$U_F$ based on a cascade of multi-controlled~$CNOT$ (i.e., multi-controlled~$X$) gates, as,  e.g., mentioned in~\cite{klappenecker2001discrete}. On a quantum computer, the multi-controlled~$CNOT$ gates are implemented by decomposing them into standard~$CNOT$ gates. This process leads to a significant increase in the number of gates and possibly the introduction of ancilla qubits~\cite{barenco1995elementary}, \cite[Ch. 4.3]{ikeAndMike}.  The proliferation of~$CNOT$ gates can be significantly reduced by switching to the circuit depicted in Figure~\ref{fig:matp_cascade}, which is based on quantum ripple-carry addition~\cite{camps2024explicit,wang2025comprehensive}. The ripple-carry incrementer uses~$CNOT$ and Toffoli gates, that is, $X$-gates with one and two controls, and ancillas for intermediate carry propagation.  

\begin{figure}[]
	\centering
	\subfloat[][Shift unitary $U_F$\label{fig:cost-matp}]{
		\includegraphics[width=0.475\textwidth]{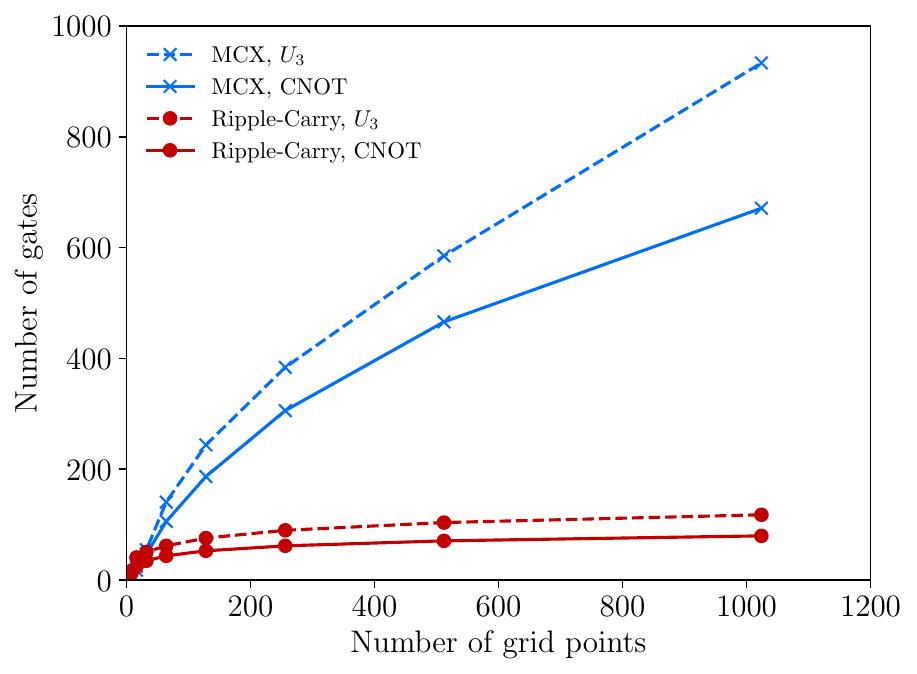}
	}
	\hfill
	\subfloat[][Reflection unitary $U_R$\label{fig:cost-matp-reflect}]{
		\includegraphics[width=0.475\textwidth]{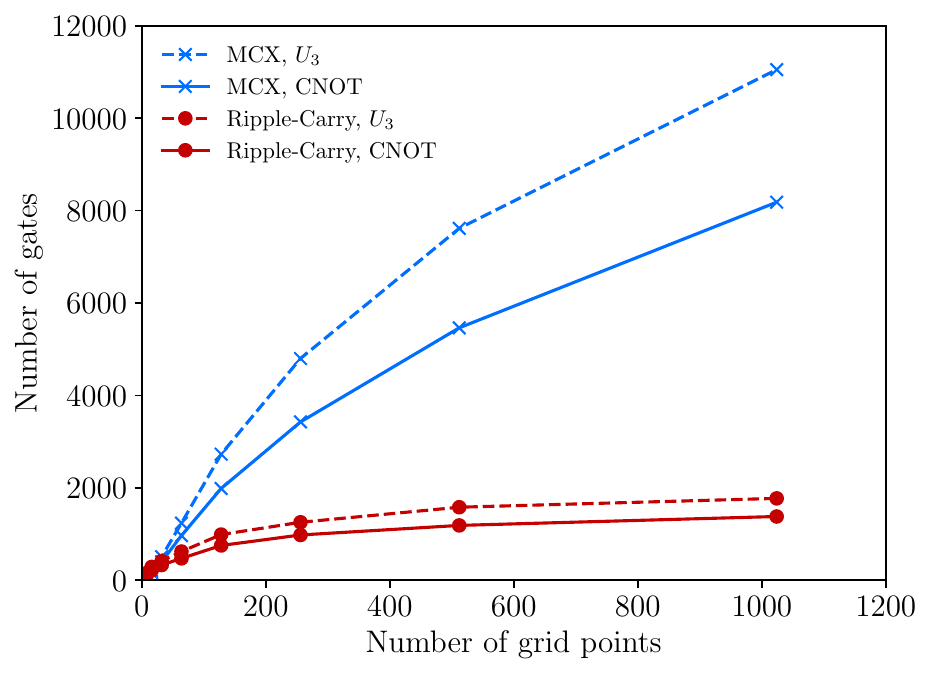}
	}	
	\caption[]{Total number of $CNOT$ and $U_3$ gates for (a) the shift unitary $U_F$ and (b) its component the reflection unitary~$U_R$. The label~$MCX$ refers to the implementation using multi-controlled~$X$ gates and the label Ripple-Carry to the implementation based on ripple-carry adder in Figures~\ref{fig:matp_mcx} and~\ref{fig:matp_cascade}, respectively.}  
	\label{fig:shift_unitary_scaling}
\end{figure}

As detailed in \ref{appx:alg_complex}, any quantum circuit can be expressed using only the two-qubit $CNOT$ and the single-qubit $ U_3$ gates. We analyse the scaling of the two different~$ U_F$ implementations by expressing the two circuits using only the $CNOT$ and $ U_3$ gates and counting the total number of gates. As evident from the scaling plots for the forward shift unitary~$U_F$ and the entire reflection unitary~$U_R$ in Figure~\ref{fig:shift_unitary_scaling}, the ripple-carry incrementer leads to significantly lower gate counts compared to the multi-controlled~$CNOT$ incrementer. In passing, we note that other implementations of the forward shift unitary are possible, for instance, using the QFT-based Draper adder~\cite{draper2000addition,schalkers2024efficient}.

%
\subsection{Solution in the spectral domain \label{sec:quantum-poly}}
%
The solution of boundary value problems in the Fourier domain requires the encoding of the solution operator given by the diagonal matrix~$\vec D^{-1}$, cf.~\eqref{eq:solution_poisson_fsp} and~\eqref{eq:solution_poisson_fractional}. Focusing on the Poisson-Dirichlet problem, we consider the encoding of the function
\begin{equation} 
	\label{eq:relabelling_3.3}
	d(k) = \begin{cases}
		\left ( \dfrac{L}{ \pi k} \right )^{2} & 1 \le k < N/2 \\[0.5em]
		0 & N/2 \le k < N \, , 
	\end{cases}
\end{equation}
representing the diagonal entries of~$\vec D^{-1}$. The function value~$d(0)$ for~$k=0$ does not affect the solution and can be arbitrarily chosen. To facilitate the encoding of~$d(k)$, we consider the function 
\begin{equation}
	\tilde d(k) =  \arcsin \left ( \frac{1}{k^2} \right ) \,.
\end{equation}
The factor~$L^2/\pi^2$ in~\eqref{eq:relabelling_3.3} can be eliminated by rescaling the forcing term by the same amount so that~$d(k)$ and~$\tilde d(k)$ are the same up to the~$\arcsin$ reparameterisation. The function~$\tilde d(k)$ is evaluated on a quantum computer via the unitary 
\begin{equation} \label{eq:poly_unitary}
	U_P (\tilde d(k))  \colon \ket k \ket 0 \mapsto  \ket k R_Y( 2 \tilde d(k))  \ket 0  =   \sqrt{1-\frac{1}{k^4}} \ket k \ket 0 + \frac{1}{k^2} \ket k \ket 1  \, ,
\end{equation}
where~$\ket k \ket 0 \equiv \ket{k_0 \dotsc k_{n-2} k_{n-1} } \ket 0 $. The rotation gate~$R_Y(2 \tilde d(k))$ is defined in terms of the matrix exponential of the Pauli $Y$ gate and is given by
\begin{equation} \label{eq:ry_gate}
	R_Y(2 \tilde d(k)  ) = e^{-i Y  \tilde d(k)} = \begin{pmatrix*}[r] \cos \tilde d(k)   & - \sin \tilde d(k) \\ \sin \tilde d(k)  &  \cos \tilde d(k)    \end{pmatrix*}   =
	\begin{pmatrix*}[c] \sqrt{1-\frac{1}{k^4}}  & - \frac{1}{k^2}  \\  \frac{1}{k^2}    &  \sqrt{1-\frac{1}{k^4}}   \end{pmatrix*}  \, .  
\end{equation}
Hence, the unitary~$U_P  (\tilde d(k))  \in  \mathbb R^{2 N \times 2 N}$ can be identified as a 2D rotation in the plane spanned by the two vectors~$\ket k \ket 0 $ and~$\ket k \ket 1 $. The last qubit in~\eqref{eq:poly_unitary}, initially in state~$\ket 0$,  is an ancilla and the amplitude~$1/k^2$ of the component~$\ket k \ket 1$ is the sought function value~$\tilde d(k)$. The amplitude of the component~$\ket k \ket 0$ is of no further interest.  The ancilla qubit is required to ensure that the output vector is of unit length, despite the scaling applied to the basis vector~$\ket k$ by~$\tilde d(k)$. It is worth noting that instead of the~$\arcsin$ parameterisation, an~$\arccos$ parameterisation is equally possible.

Next, we approximate~$\tilde d(k)$ by a polynomial
\begin{equation} \label{eq:poly_approx}
	\tilde  d(k)  \approx \sum_{j=0}^p  \alpha_j k^j \, , 
\end{equation}
where~$p$ is the chosen degree. The coefficients~$\alpha_j$ are determined on a classical computer using standard polynomial approximation techniques. 
With a slight abuse of notation, we use the same symbol for~$\tilde d(k)$ and its polynomial approximation in what follows. Moreover, we assume that the interval~$1 \le k \le N/2$ is approximated with only one polynomial. In our computations, we also consider the decomposition of the interval into several sub-intervals; see~\ref{appx:bivariate_piecewise} for the treatment of multiple intervals and multivariate polynomials.  Expressing the argument~$k$ of the polynomial~\eqref{eq:poly_approx} as a binary and using the multinomial theorem yields
\begin{align}
\begin{split}
	\tilde  d(k) &=  \sum_{j=0}^p  \alpha_j \left (k_0 2^{n-1} + \cdots + k_{n-2} 2^1 + k_{n-1} 2^0   \right )^j  = \sum_{j=0}^p  \alpha_j  \left (\sum_{s=0}^{n-1} k_s 2^{n-1-s}  \right )^j   \\ &=   \sum_{j=0}^p \alpha_j \left ( \sum_{r_0+r_1+\cdots+r_{n-1} =j}  \frac{j}{r_0! r_1! \cdots r_{n-1}!} \prod_{s=0}^{n-1} \left (  k_s 2^{n-1-s}\right )^{r_s} \right ) \, .
\end{split}
\end{align}
Note that~$k_s = (k_s)^{r_s}$ when~$r_s \neq 0$ (since~$k_s \in \{ 0,1\}$). Hence, we can write
\begin{equation}
	\tilde  d(k) =  \alpha_0 + \sum_{j=1}^p \alpha_j \left ( \sum_{r_0+r_1+\cdots+r_{n-1} =j}  \frac{j}{r_0! r_1! \cdots r_{n-1}!} \prod_{s=0}^{n-1} k_s 2^{r_s(n-1-s)}  \right )  \, .
\end{equation}
Introducing this expression in~\eqref{eq:poly_unitary} and~\eqref{eq:ry_gate}, the summation reduces to a multiplication owing to the commutativity of 2D rotations. The rotation angles depend on~$k_s$ representing the state of the qubit~$s$, hence leading to controlled rotations.

\begin{figure}[]
	\centering
\begin{quantikz}[row sep=0.5cm, column sep=0.5cm]
	\lstick{$\ket{k_0}$}   & \qw &\qw  
	\gategroup[3, steps=1,style={rounded
		corners,fill=csmlPink, opacity=0.5, line width =0.1pt, inner xsep=0.pt},background, label style={label
		position=below,anchor=north,yshift=-0.45cm}]{$U_{P_0}$}  
	&  \ctrl{2}
		\gategroup[3, steps=1,style={rounded
		corners,fill=csmlPink, opacity=0.5, line width =0.1pt,inner xsep=0.pt},background, label style={label
		position=below,anchor=north,yshift=-0.45cm}]{$U_{P_1}$}  
	 & \qw
	 	\gategroup[3, steps=1,style={rounded
		corners,fill=csmlPink, opacity=0.5, line width =0.1pt, inner xsep=0.pt},background, label style={label
		position=below,anchor=north,yshift=-0.45cm}]{$U_{P_2}$}  
	  &   \ctrl{2} 
	  	\gategroup[3, steps=1,style={rounded
		corners,fill=csmlPink, opacity=0.5, line width =0.1pt, inner xsep=0.pt},background, label style={label
		position=below,anchor=north,yshift=-0.45cm}]{$U_{P_3}$}  
	  & \qw & \qw \\
	\lstick{$\ket{k_1}$ } &\qw & \qw & \qw  & \ctrl{1}   & \ctrl{1} & \qw & \qw  \\
	\lstick{$\ket{0}$}  &\qw  & \gate[label style={yshift=0.cm}, style={fill=csmlGrey9, opacity=0.8}]{R_Y(2\alpha_0)} &  \gate[label style={yshift=0.cm}, style={fill=csmlGrey9, opacity=0.8}]{R_Y(4 \alpha_{1}+8\alpha_{2})}  &  \gate[label style={yshift=0.cm}, style={fill=csmlGrey9, opacity=0.8}]{R_Y(2 \alpha_{1}+2\alpha_{2})}  & \gate[label style={yshift=0.cm}, style={fill=csmlGrey9, opacity=0.8}]{R_Y(8 \alpha_2)} 
	& \qw& \qw 
\end{quantikz}
	\caption{Quantum circuit for the unitary~$U_P$ evaluating the function~$\tilde d(k) = \alpha_0 + \alpha_1 k + \alpha_2 k^2$ for~$k \in \{0, \, 1, \, 2, \, 3 \}$, where~$k = k_0 2^1 + k_1 2^0$ and~$k_0, k_1 \in \{ 0, \, 1\}$. The circuit is composed of four unitaries with each unitary implementing a 2D rotation in a plane defined by the vectors~$\ket{k} \ket 0$ and~$\ket{k} \ket 1$, equivalently~$\ket{k_0 k_1} \ket 0$ and~$\ket{k_0 k_1} \ket 1$.  \label{fig:polycirc_quad}}
\end{figure}
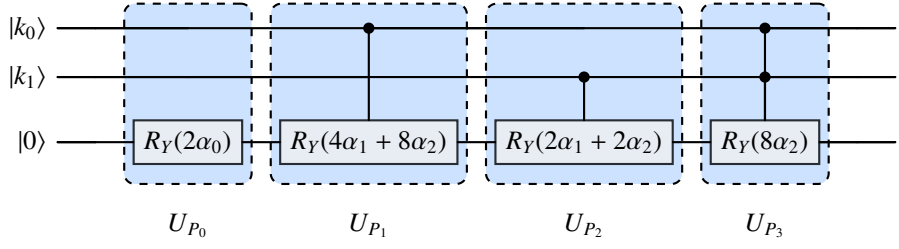

As an example, Figure~\ref{fig:polycirc_quad}  shows the circuit for evaluating the quadratic polynomial
\begin{equation}
	\tilde d(k) = \alpha_0 + \alpha_1 k + \alpha_2 k^2 \, 
\end{equation}
for~$k \in \{ 0, \, 1, \, 2, \, 3 \}$, i.e.,~$k = k_0 2^1 + k_1 2^0$ with~$k_0, k_1 \in \{0, 1 \}$. This function expands to
\begin{equation}
	\tilde d(k_0,k_1) = \alpha_0 + k_0 (2 \alpha_1 + 4 \alpha_2) + k_1 (\alpha_1 + \alpha_2) + k_0k_1(4 \alpha_2) \, .
\end{equation}	
Consequently, the circuit can be implemented using the controlled rotations~$R_Y(2 \alpha_0)$, $R_Y(4 \alpha_1 + 8 \alpha_2 )$, $R_Y(2\alpha_1 + 2 \alpha_2)$ and $R_Y(8 \alpha_2)$, or the 
multiplicative decomposition 
\begin{equation}
	U_P (\tilde d(k_0,k_1)) = U_{P_3} (8 \alpha_2) U_{P_2} (2\alpha_1 + 2 \alpha_2) U_{P_1}(4 \alpha_1 + 8 \alpha_2)  U_{P_0} (2 \alpha_0)  \, ,
\end{equation}
with each unitary implementing a 2D rotation in a specific plane given by the vectors~$\ket {k_0 k_1}\ket 0$ and~$\ket {k_0 k_1} \ket 1$. For instance, the first two unitaries are given by 
\begin{equation}
	U_{P_0} (2 \alpha_0)  =  I \otimes I \otimes R_{Y}(2 \alpha_0) = \begin{pmatrix} R_{Y}(2 \alpha_0)  & & & \\ &  R_{Y}(2 \alpha_0) & & \\   &  &  R_{Y}(2 \alpha_0)  & \\ &  & & R_{Y}(2 \alpha_0)  \end{pmatrix} \, ,
\end{equation}
\begin{equation}
	U_{P_1}  (4 \alpha_1 + 8 \alpha_2 ) =  \ket 0 \bra 0   \otimes I \otimes  I  + \ket 1 \bra 1  \otimes I \otimes   R_{Y}(4 \alpha_1 + 8 \alpha_2) = \begin{pmatrix} I \quad  & & & \\ & I \quad  & & \\   &  &  R_{Y}(4 \alpha_1 + 8 \alpha_2)  & \\ &  & & R_{Y}(4 \alpha_1 + 8 \alpha_2)   \end{pmatrix} \, .
\end{equation}
%
\subsection{Quantum spectral solver \label{sec:quantum-solver}}
%
We are now ready to introduce the quantum circuit for solving Dirichlet boundary value problems by combining state preparation, QFT, reflection, and polynomial encoding. For state preparation, any of the numerous state preparation algorithms can be used. The other necessary algorithms have all been discussed in this section.  We first introduce the circuit for one-dimensional problems and then briefly extend it to two-dimensional problems.  

\begin{figure}[]
	\centering
	\input{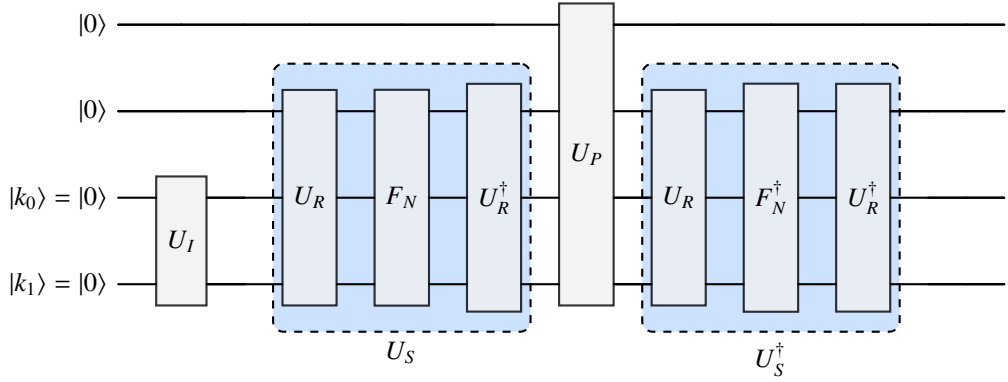}
	\caption{Quantum circuit for solving a one-dimensional Dirichlet boundary value problem. The extended domain~$\Omega_{\text E}$ is discretised with~\mbox{$N=2^3=8$} grid points so that the physical domain~$\Omega$ is discretised with $4$ grid points. The unitary~$U_I$ encodes the source~$\ket f$ corresponding to the physical domain, and the output of the circuit is the solution~$\ket u$ of the boundary value problem.}
	\label{fig:quantum_solver_1d}
\end{figure}

The one-dimensional circuit in Figure~\ref{fig:quantum_solver_1d} uses~$n=3$ qubits to represent the field vectors on the extended domain and one ancilla qubit for polynomial encoding. The vectors corresponding to the extended domain have~$N=2^n$ components, and the ones corresponding to the problem domain have~$N/2$ components and~$n-1$ qubits. The unitary~$U_S =U_R^\dagger F_N U_R $ implements the quantum sine transform and is composed of the unitary matrices~$U_R$ implementing the antisymmetric reflection and~$F_N$ implementing the quantum Fourier transform for an~$N$-dimensional vector.  
All the qubits in the circuit are initialised in state~$\ket 0$. Using the state preparation unitary~$U_I$, the prescribed source vector is encoded into the amplitudes of the state vector 
\begin{equation}
	\label{eq:fstate}
	I^{\otimes 2} \otimes U_I \colon \quad \ket {0}^{\otimes 2}  \ket {0}^{\otimes 2}   = \sum_{k=0}^{N/2-1}  f_k   \ket {0}^{\otimes 2}  \ket k \, , 
\end{equation}
where~$\ket k \equiv \ket{k_0 k_1} \equiv \ket {k_0} \ket {k_1}$, with~$k \in \{0,\, 1, \, 2, \, 3\}$ and~$k_0, k_1 \in \{ 0,  \, 1 \}$. The sequence of subsequent mappings is as follows
\begin{subequations}
\begin{align}
	I \otimes U_S \colon&  \quad    \sum_{k=0}^{N/2-1} f_k   \ket {0}^{\otimes 2}  \ket k  \mapsto  \sum_{k=0}^{N/2-1} \hat f_k   \ket {0}^{\otimes 2}  \ket k  \, ,\\
	U_P \colon &  \quad  \sum_{k=0}^{N/2-1} \hat f_k   \ket {0}^{\otimes 2}  \ket k \mapsto  \sum_{k=0}^{N/2-1} \hat f_k   \cos \tilde d(k) \ket 0 \ket 0 \ket k   +  \sum_{k=0}^{N/2-1} \hat f_k   \sin \tilde d(k) \ket 1 \ket 0 \ket k  \, ,  \\
	I \otimes U_S^\dagger \colon & \quad  \sum_{k=0}^{N/2-1} \hat f_k   \cos \tilde d(k) \ket 0 \ket 0 \ket k   +  \sum_{k=0}^{N/2-1} \hat f_k   \sin \tilde d(k) \ket 1 \ket 0 \ket k  \mapsto   \sum_{k=0}^{N/2-1} u_k   \ket 1 \ket 0 \ket k + \text{Junk} \, .
\end{align}
\end{subequations}
After the circuit is executed, the solution of the boundary value problem corresponds to the amplitudes of the states with the first (top) qubit in state~$\ket 1$. The irrelevant terms with the first qubit in state~$\ket 0$ are denoted as ``Junk''. The second qubit is only used within the implementation of~$U_S$ and will be in state~$\ket 0$ outside of~$U_S$ (on a noiseless quantum computer).  Although we assumed in this section that the solution in the Fourier space is encoded using a single polynomial, in practice, we use a piecewise polynomial approximant. This affects only the design of the unitary~$U_P$ and requires ancilla qubits beyond the ones present in Figure~\ref{fig:quantum_solver_1d}.

\begin{figure}[]
	\centering
	\input{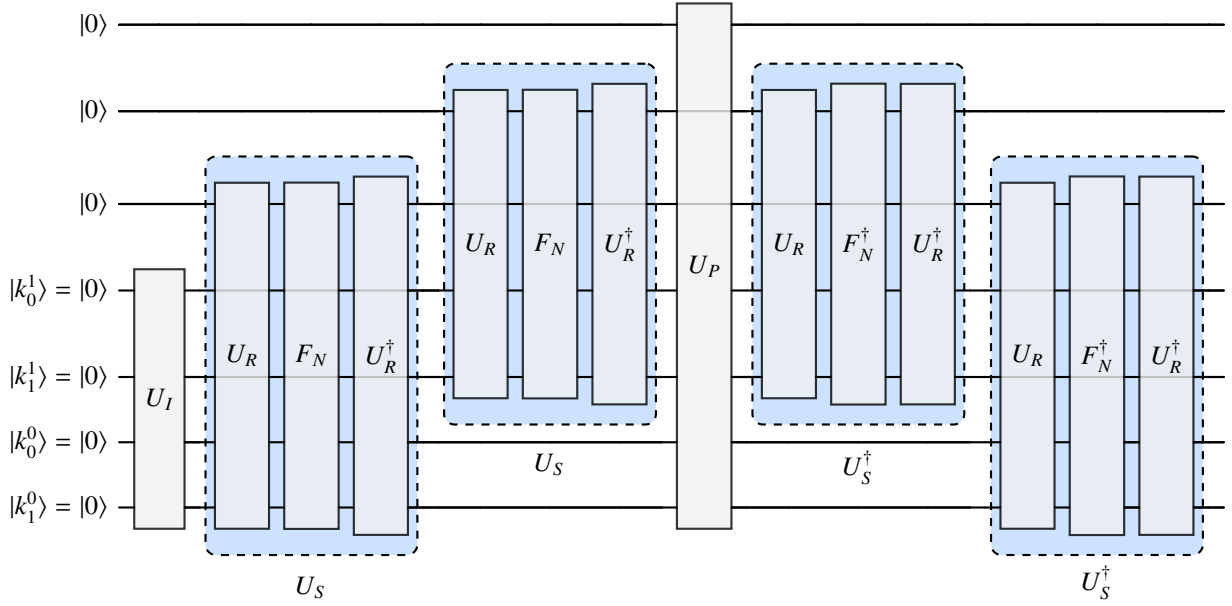}
	\caption{Quantum circuit for solving a two-dimensional Dirichlet boundary value problem. The extended domain~$\Omega_{\text E}$ is discretised with~\mbox{$N \times N =2^3 \times 2^3 =8 \times 8$} grid points so that the physical domain~$\Omega$ is discretised with $4 \times 4$ grid points. The unitary~$U_I$ encodes the source~$\ket f$ corresponding to the physical domain, and the output of the circuit is the solution~$\ket u$ of the boundary value problem.  \label{fig:solver2d}}
\end{figure}

The circuit diagram for two-dimensional Dirichlet boundary value problems is shown in Figure~\ref{fig:solver2d}. The extended domain is discretised with~$N\times N = 8 \times 8$ grid points. That is, the field vectors of the extended domain are encoded as the state vector of~$n\times n = 3 \times 3$ qubits and the ones of the problem domain as the state vector of~$n-1 \times n-1$ qubits. The grid points are labelled with the multi-index 
\begin{equation}
	\ket {\vec k} \equiv \ket {k^0} \ket{k^1}  \equiv  \ket {k_0^0} \ket {k_1^0} \ket {k_0^1} \ket {k_1^1} \, ,  \quad k^0, k^1 \in \{ 0, \,  1,  \, 2, \, 3 \} \, ,  \quad  k_0^0, k_1^0, k_0^1, k_1^1 \in \{ 0, \, 1 \} \,. 
\end{equation}
The two-dimensional solver is constructed mostly as the tensor product of two one-dimensional solvers.  In  Figure~\ref{fig:solver2d}, the unitary~$U_S$ is exactly the same as for one-dimensional problems; it is applied simultaneously in the coordinate directions~$ k^0$ and~$ k^1$. Although not evident from the circuit diagram, the qubits for the two coordinate directions~$ k^0$ and~$ k^1$ are independent from each other so that~$U_S$ is applied simultaneously. The primary difference between the one- and two-dimensional problems concerns the solution of the boundary value problem in the spectral domain. Hence, the unitary~$U_P$ is different from the one-dimensional case and depends on the specific boundary value problem considered. 

\section{Examples}
\label{sec:example}
%
We present a series of numerical experiments to assess the convergence and computational complexity of the proposed quantum approach. All circuits are implemented in Qiskit and executed on a noiseless simulator~\cite{qiskit2024}. Convergence is studied by comparing the $L_2$-norm error of the quantum solution with that of the classical spectral solution. We assess computational complexity by expressing the circuits in terms of only two-qubit~$CNOT$ gates and single-qubit rotation gates~$U_3$. The gate set $\{ CNOT, \, U_3\}$ is universal, meaning that any quantum circuit can be expressed using these two gates. We express the circuits with the~$\{ CNOT, \, U_3\}$ gate set by compiling them either in Qiskit or with TKET~\cite{qiskit2024,tket2020}, see~\ref{appx:alg_complex} for details.  TKET usually yields significantly lower gate counts by applying different algorithms for circuit optimisation. We first study the quantum solution of one-dimensional boundary value problems, specifically the Poisson-Dirichlet problem with both homogeneous and inhomogeneous boundary conditions, as well as a fractional differential equation. Subsequently, we study the quantum solution of these problems on two-dimensional domains. 
%
\subsection{One-dimensional boundary value problems \label{sec:onedim}}

\subsubsection{Homogeneous Poisson-Dirichlet problem \label{sec:oneD-poisson-homog}}
\begin{figure}[]
	\centering
	\subfloat[][Forcing \label{fig:1d-poisson-force}]{
		\includegraphics[width=0.475\textwidth]{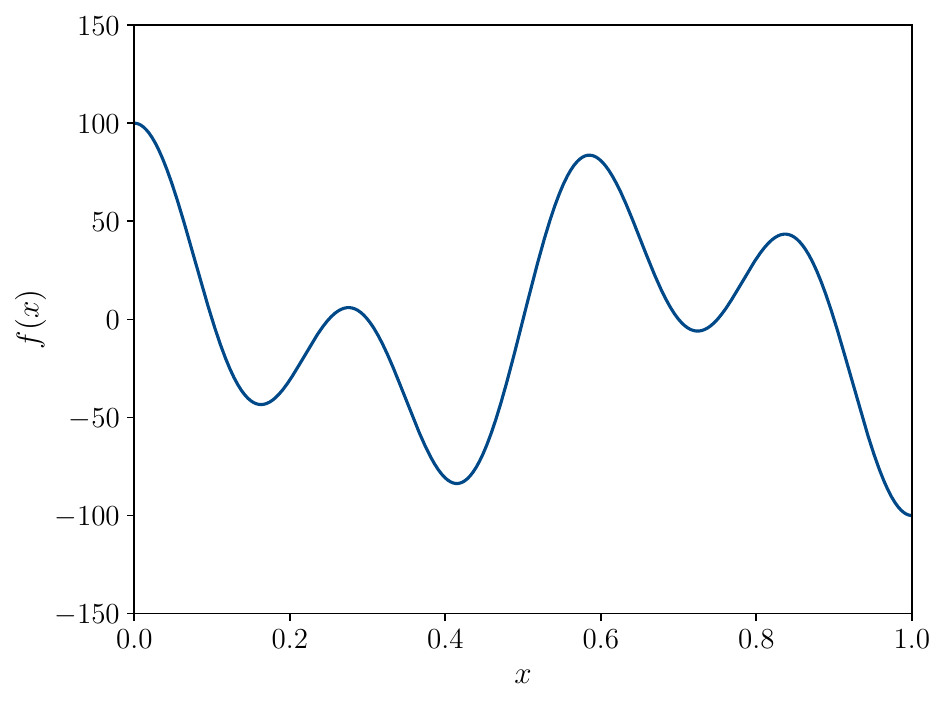}
	}
	\hfill
	\subfloat[][Solution \label{fig:1d-poisson-soln}]{
		\includegraphics[width=0.475\textwidth]{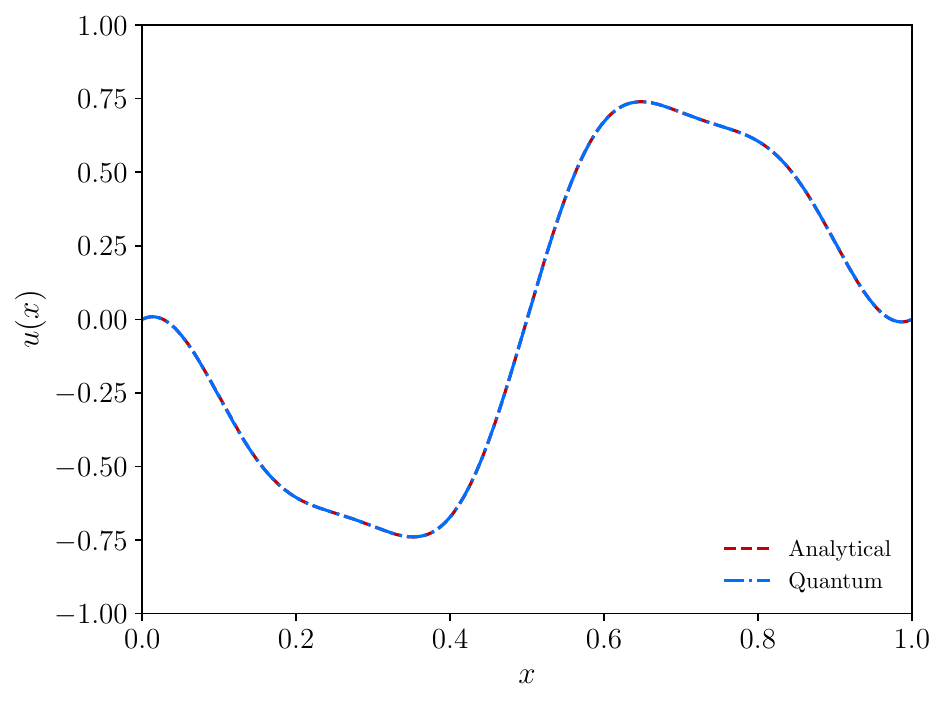}
	}
	\caption[]{One-dimensional homogeneous Poisson-Dirichlet problem. \label{fig:1d-poisson-f-u}}
\end{figure}

%
We seek in the domain~$\Omega = (0, \,1 )$ the solution of the  Poisson-Dirichlet boundary value problem with the boundary conditions $u(0) = u(1) = 0$ and the source field 
\begin{equation}
	\label{eq:force-1d-poi}
	f(x) =   100 \, \cos \left ( 2 \pi x \right) \, \cos \left( 5 \pi x \right) \, .
\end{equation}
Its analytical solution is given by 
 \begin{equation}
 	\label{eq:anlyt-1d-poi}
 	u(x) = \frac{ 100 \left(-58 + 116 x + 49 \cos \left( 3 \pi x \right) + 9 \cos \left( 7 \pi x\right) \right)}{882 \pi^2} \, ; 
 \end{equation}
see Figure~\ref{fig:1d-poisson-f-u}. The domain $\Omega$ is discretised with $N$ uniformly distributed grid points with a spacing of $h = 1/N$. The problem is then solved with the quantum circuit introduced in Section~\ref{sec:quantum-solver} by applying a sequence of unitary transformations. The input quantum state vector $\ket{f}$ is obtained by evaluating $f(x)$ at the grid points and encoding it using a state preparation algorithm. The state vector~$\ket f$ must have a unit length, which is achieved by scaling~$f(x)$. The computed solution is rescaled to account for the initial scaling of~$f(x)$. The solution of the Poisson-Dirichlet problem reduces to algebraic equations in the Fourier space. We encode the solution operator~$\vec D^{-1}$ introduced in~\eqref{eq:solution_poisson_fsp} using piecewise polynomials of degrees~$p=3$ and~$ p=4$ as described in Section~\ref{sec:quantum-poly}. To improve the approximation of~$\vec D^{-1}$, we choose smaller segment sizes closer to the origin. The segment size choice is optimised by trial and error and has not been further formalised in this work.

\begin{figure}[]
	\centering
	\subfloat[][$L^2$-norm error \label{fig:1d-poisson-conv-H}]{
		\includegraphics[width=0.475\textwidth]{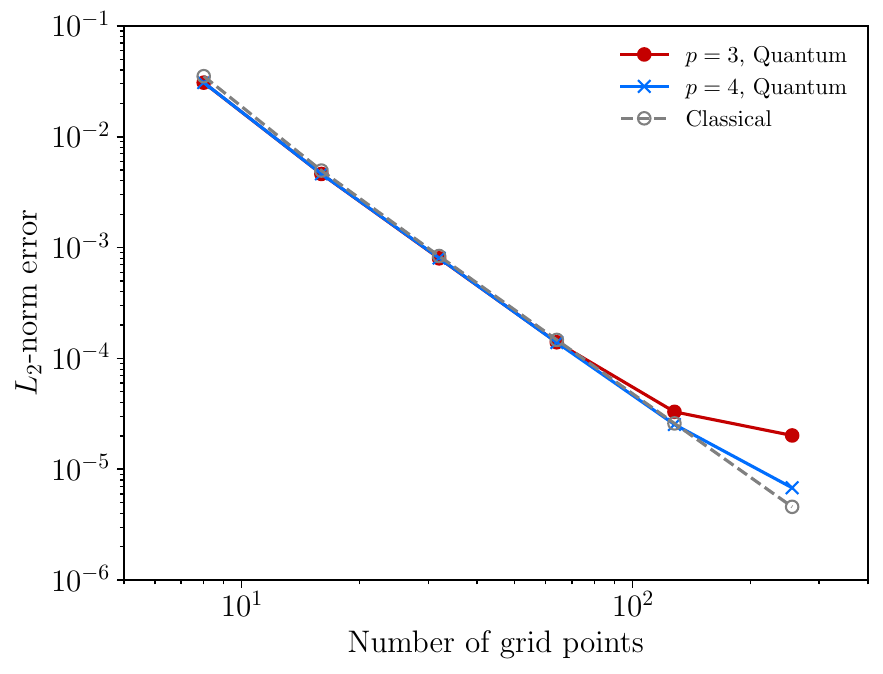}
	}
	\hfill
	\subfloat[][Total number of~$CNOT$ and~$U_3$ gates \label{fig:1d-poisson-cost-H}]{
		\includegraphics[width=0.475\textwidth]{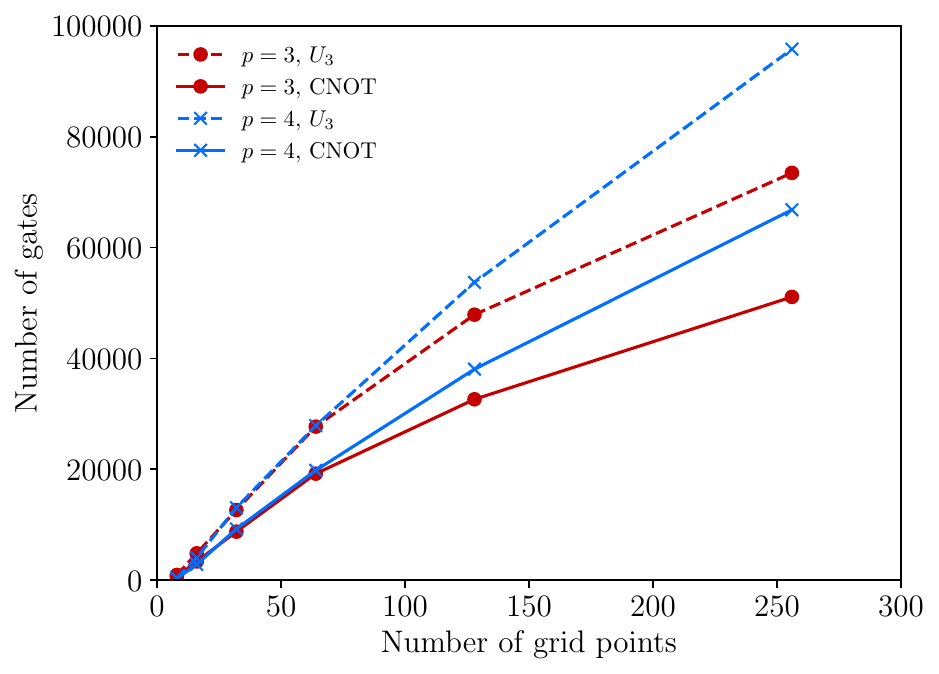}
	}
	\caption[]{One-dimensional homogeneous Poisson-Dirichlet problem. \label{fig:1d-poisson-h}}
\end{figure}

The convergence of the classical and quantum spectral approximation is shown in  Figure~\ref{fig:1d-poisson-conv-H}. The classical approximation exhibits uniform convergence, whereas the convergence of the quantum approximation exhibits a flattening when refining the grid and depends on the degree~$p$ in the polynomial approximation of the solution operator~$\vec D^{-1}$. However, the accuracy of the quantum solution for finer meshes can be systematically increased by choosing either a higher~$p$ or smaller segment sizes in piecewise polynomial approximation of~$\vec D^{-1}$. As visible in Figure~\ref{fig:1d-poisson-cost-H}, $p$  has an effect on the number of $CNOT$ and~$U_3$ gates in the quantum circuit. Higher polynomial degrees and smaller segment sizes lead to an increase in gate counts. Irrespective of their total number, the number of gates always depends polylogarithmically on the number of grid points. 

\subsubsection{Inhomogeneous Poisson-Dirichlet problem \label{sec:oneD-poisson-inhomog}}
\begin{figure}[]
	\centering
	\subfloat[][Solution $u(x)$ and auxilliary function $g(x)$ \label{fig:1d-popisson-nh-a}]{
	\includegraphics[width=0.475\textwidth]{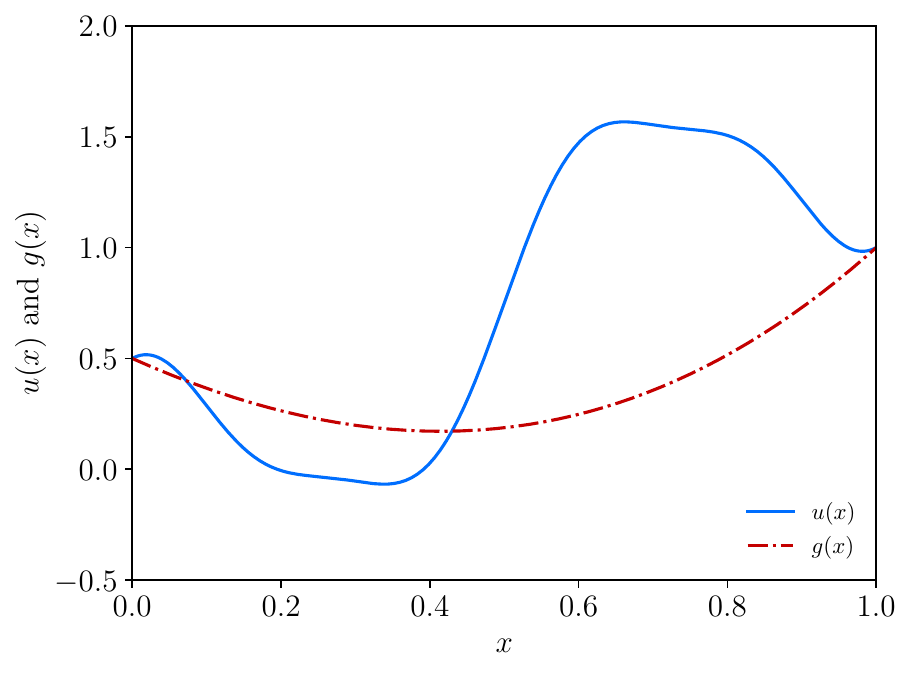}
	}
	\hfill
	\subfloat[][$L^2$-norm error \label{fig:1d-poisson-conv-NH} \label{fig:1d-popisson-nh-b}]{
		\includegraphics[width=0.475\textwidth]{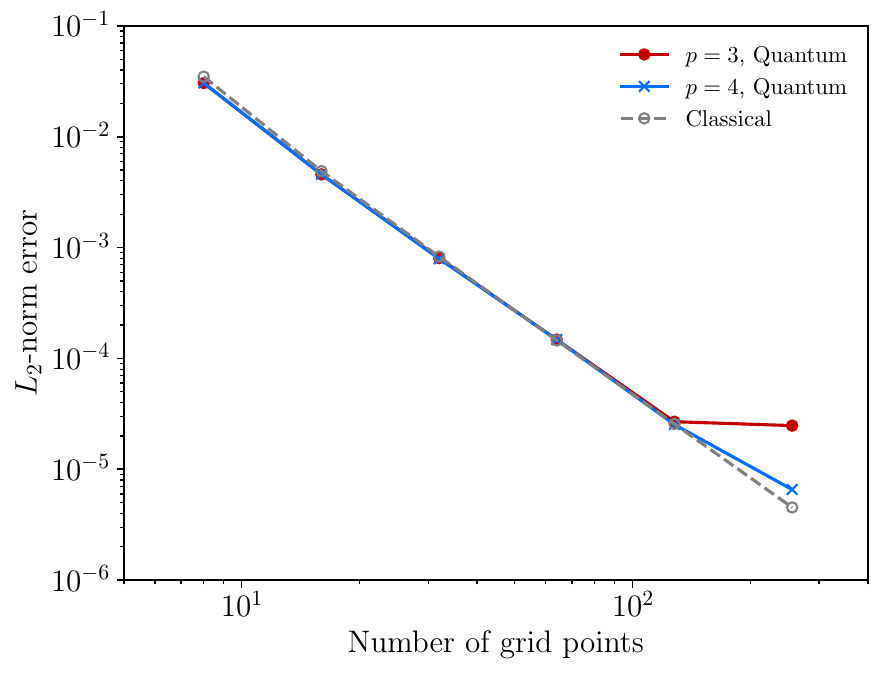}
	}
	\caption[]{One-dimensional inhomogeneous Poisson-Dirichlet problem. \label{fig:1d-popisson-nh}}
\end{figure}

%
Building on the previous example,  we consider the inhomogeneous Poisson-Dirichlet problem on the domain $\Omega = (0, 1)$.  The boundary conditions are now~$u(0) = 0.5$ and $u(1) = 1$, and the source field~$f(x)$ is the same as in~\eqref{eq:force-1d-poi}.  The analytical solution becomes
 \begin{equation}
	\label{eq:anlyt-1d-poi-nh}
	u(x) = 0.5(x+1)+  \frac{ 100 \left(-58 + 116 x + 49 \cos \left( 3 \pi x \right) + 9 \cos \left( 7 \pi x\right) \right)}{882 \pi^2} \, .
\end{equation}
To solve this problem, we consider the homogeneous Poisson-Dirichlet problem 
\begin{equation} \label{eq:homog_dirichlet_inhomog}
\begin{aligned}
	& - \frac{\D^2 v(x)}{\D x^2} = f(x) + \frac{\D^2 g(x)}{\D x^2}  \, , \quad  \forall x \in \Omega \, , \\
	& v(0) = 0 \, , \quad  v(1) = 0  \, . 
\end{aligned}
\end{equation}
The auxiliary function is chosen as 
\begin{equation}
	g(x) =  - 0.5(x-1)^3 + x^3  \, ,
\end{equation}
 and satisfies the Dirichlet boundary conditions and is sufficiently smooth. After computing~$v(x)$ the solution of the original problem is recovered using the relationship~$u(x) = v(x)+g(x)$, see Figure~\ref{fig:1d-popisson-nh-a}.

The homogeneous Poisson-Dirichlet problem for~$v(x)$  is solved using the classical or quantum spectral methods following the procedure in the previous example. In the quantum approach, we first compute the modified source term classically and then classically recover the solution from the quantum solution of the homogeneous problem. The convergence of the classical and quantum spectral approximations is shown in Figure~\ref{fig:1d-popisson-nh-b}. 

%
\subsubsection{Fractional stochastic differential equation \label{sec:oneD-fractional}}
%
We next consider a problem motivated by the representation of random fields as the solution of fractional stochastic differential equations. For convenience, we restate the fractional differential equation~\eqref{eq:fractional_ode} as 
\begin{equation}  \label{eq:fractional_ode_ex}
		\left ( \kappa^2  - \frac{\D^2 }{\D x^2} \right )^\beta  u ( x) = \frac{1}{\tau} f (x) \, ,  \quad \forall  x \in (0, 1) \, . 
\end{equation} 
Here, the source~$f(x)$ is a Gaussian white-noise field
\begin{equation}
	f(x) \sim \mathcal{GP}\left( 0, \delta (x - x') \right) \, , 
\end{equation}
with the Dirac delta covariance function~$\delta(x-x')$. As is well known, see e.g.~\cite{koh2023stochastic,lindgren2011explicit}, on an unbounded domain the solution~$u(x)$ is a Gaussian process 
\begin{equation}
	u(x) \sim \mathcal{GP}\left( 0, c_u(x,x') \right) \, , 
\end{equation}
where~$c_u(x, x')$ is the Mat\'ern covariance function which has a closed-form; see~\ref{appx:ode_matern}. The length-scale, smoothness and variance of the random solution~$u(x)$ are governed by the parameters~$\kappa$,~$\beta$ and~$\tau$, respectively. 

\begin{figure}[]
	\centering
	\subfloat[][Samples from a Gaussian process \label{fig:random-field}]{
			\includegraphics[width=0.475\textwidth]{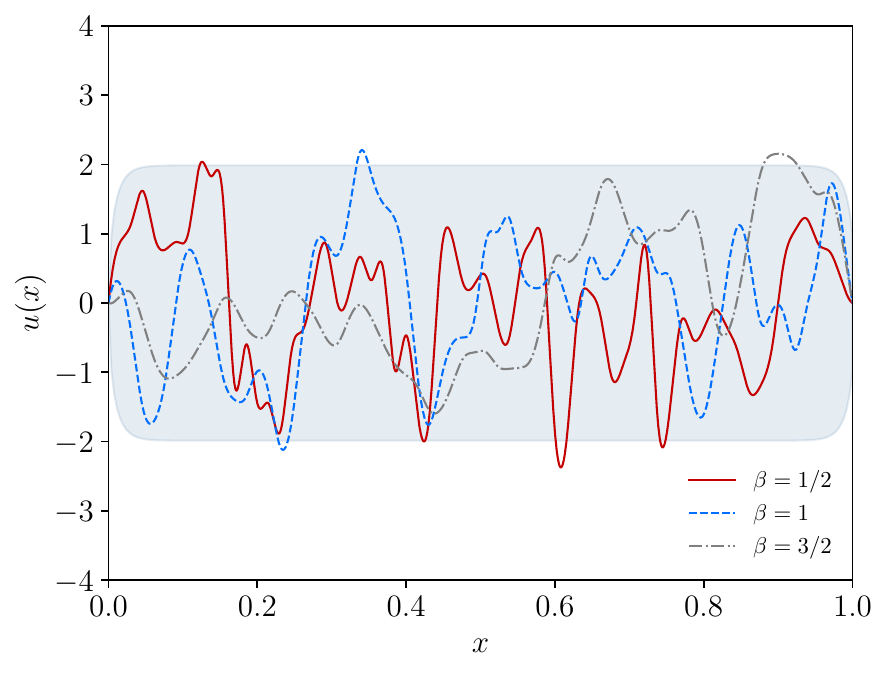}
		}
	\hfill 
		\subfloat[][Covariance functions \label{fig:1d-matern}]{
			\includegraphics[width=0.475\textwidth]{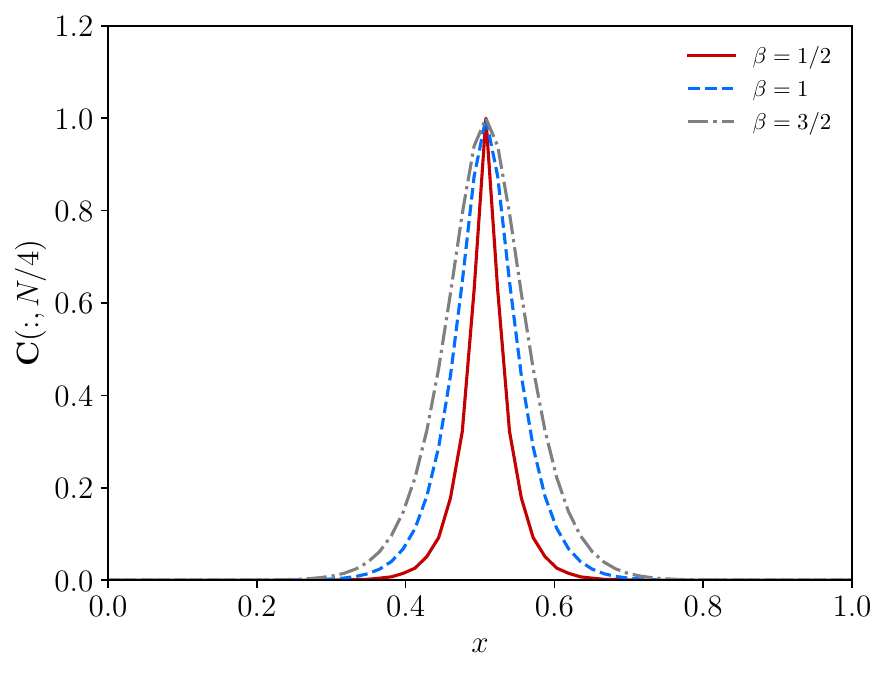}
		}
	\caption[]{One-dimensional fractional stochastic differential equation. \label{fig:1d-stoc-problem-def}}
\end{figure}

The fractional stochastic differential equation~$\eqref{eq:fractional_ode_ex}$ provides a means of sampling from the Gaussian process~$u(x)$. Specifically, a sample of~$u(x)$ is obtained by first drawing a sample from the Gaussian white-noise field~$f(x)$ and then solving~$\eqref{eq:fractional_ode_ex}$.  The discretisation of the Gaussian white-noise field yields the random vector
\begin{equation}
	\vec f = \frac{1}{N} \begin{pmatrix}  0  &  f(x_1) &  \dotsc &  f(x_{N-1}) \end{pmatrix}^\trans \,  ,
\end{equation}
where each component~$f_k \sim \mathcal{N}(0, 1/\tau^2)$ is a univariate Gaussian random variable with variance~$1/\tau^2$, and the entire vector is scaled by dividing by the number of grid points~$N$.  The respective solution vector~$\vec u$ is obtained from the quantum equivalent of~\eqref{eq:solution_fractional}, restated for convenience,  
\begin{equation*} 
	\vec u  =   \vec R^\trans \vec F^\dagger \vec R {\vec D}^{-1} \vec R^\trans \vec F \vec R \vec f \, . 
\end{equation*} 
In Figure~\ref{fig:random-field}, three sampled functions for~$\kappa=40$, \mbox{$\tau=4.279 \cdot 10^{-5}$}, and three different~$\beta \in \{1/2, \, 1, \, 3/2 \}$  are depicted.  As expected, the sample becomes progressively smoother with increasing~$\beta$. 

\begin{figure}[]
	\centering
	\subfloat[][$L^2$-norm error for $p = 3$ \label{fig:conv-p3}]{
		\includegraphics[width=0.475\textwidth]{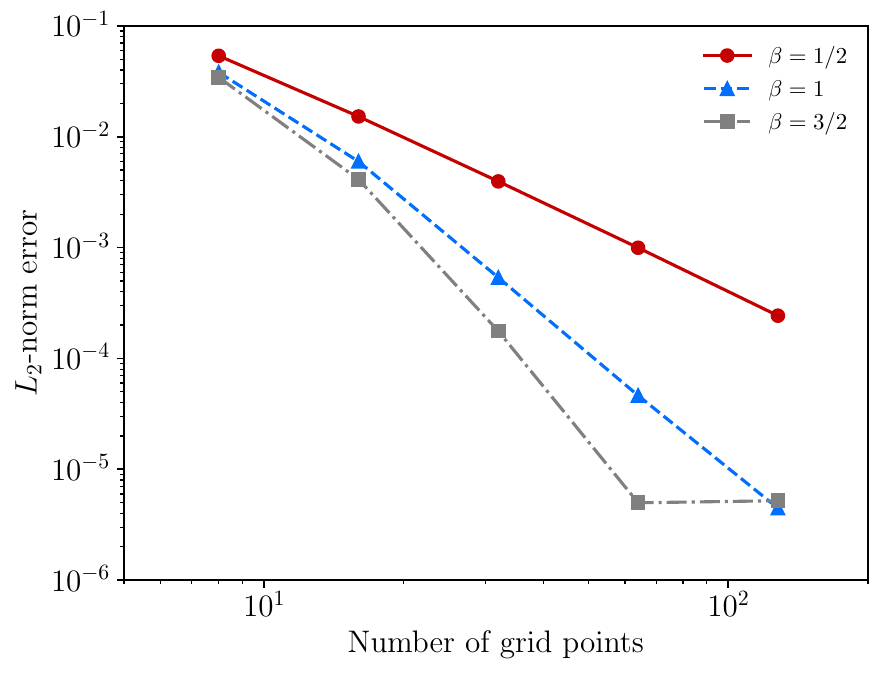}
	}
	\hfill
	\subfloat[][$L^2$-norm error for $p = 4$ \label{fig:conv-p4}]{
		\includegraphics[width=0.475\textwidth]{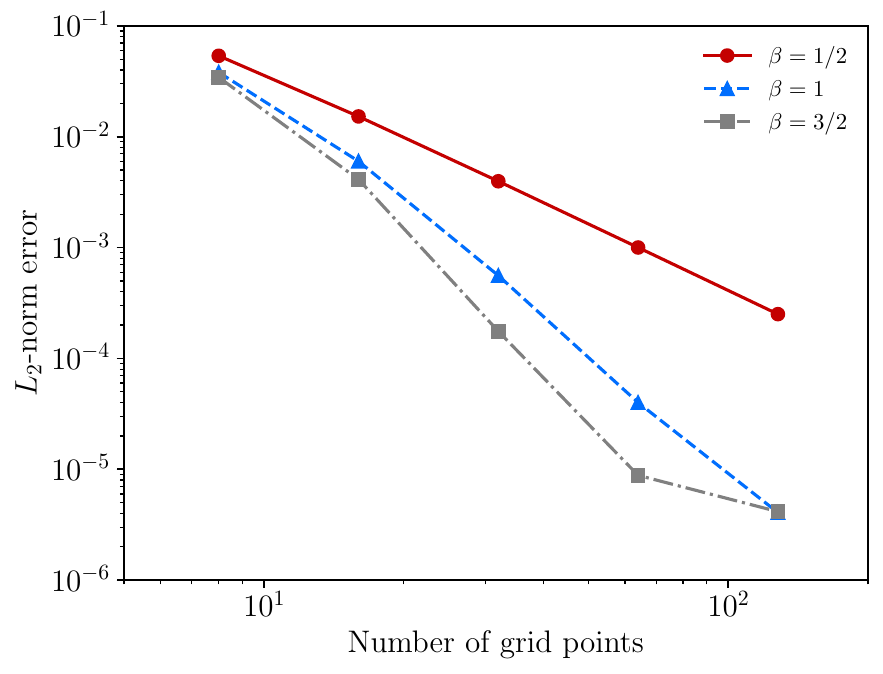}
	}
	\caption[]{One-dimensional fractional stochastic differential equation. Convergence of the $L^2$-norm error in approximating~$c_u(x, x'=1/2)$. \label{fig:1d-stoc-conv}}
\end{figure}

\begin{figure}[]
	\centering
	\includegraphics[width=0.475\textwidth]{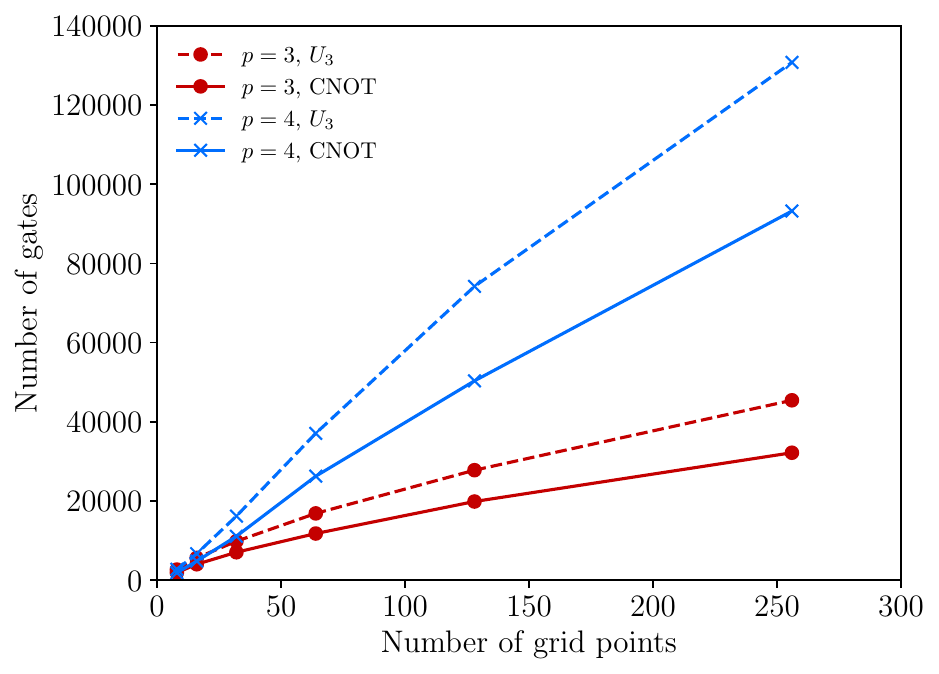}
	\caption{One-dimensional fractional stochastic differential equation. The total number of~$CNOT$ and~$U_3$ gates. \label{fig:1d-stoc-cost}}
\end{figure}

It can be shown that the discrete solution vector~$\vec u$ has the multivariate Gaussian probability density 
\begin{equation}
	\vec u \sim  \mathcal{N} (\vec 0, \, \vec C) = \mathcal{ N} \left ( \vec 0, \, \frac{1}{N}  \vec R^\trans \vec F^\dagger \vec R {\vec D}^{-2} \vec R^\trans \vec F \vec R \right ) \, ,
\end{equation}
where $\vec C$ is the covariance matrix. The entries of~$\vec C$ express the correlation between the components of the solution vector and correspond to the Mat\'ern covariance function~$c_u(x, \, x')$ of the continuous problem. In Figure~\ref{fig:1d-matern}, the central row of the covariance matrix~$\vec C(:, N/4)$, corresponding to the function~$c_u(x, \, x' =1/2)$, for the three different considered smoothness parameters is plotted. The $L_2$-norm of the difference between~$c_u(x, \, x' =1/2)$ and its approximation by~$\vec C$ for different numbers of grid points is shown in Figure~\ref{fig:1d-stoc-conv}. The convergence depends on the smoothness parameter~$\beta$, with higher convergence rates for larger~$\beta$. In the quantum approximation of the operator~$\vec D^{-2}$, we consider piecewise polynomials of degrees~$p=3$ and~$p=4$. The flattening in the convergence curves can again be attributed to the polynomial approximation error. The total number of the $CNOT$ and~$U_3$ gates in dependence of the total number of grid points~$N$ is plotted in Figure~\ref{fig:1d-stoc-cost}, confirming again their polylogarithmic scaling. 

%
\subsection{Two-dimensional boundary value problems \label{sec:twodim}}

\subsubsection{Inhomogeneous Poisson-Dirichlet problem \label{sec:twoD-Poisson}}
%
We proceed now to study the quantum spectral solution of two-dimensional boundary value problems. To begin with, we seek on the domain $\Omega = (0, \, L)^2 \subset \mathbb R^2$ the solution of the Poisson-Dirichlet problem with the source field
\begin{equation}
	f(\vec x) = 13 \pi^2  \sin(2 \pi x_0) \sin(3 \pi x_1) + 17 \pi^2 \sin(\pi x_0) \sin(4 \pi x_1) - 9 x_0 - 15 x_1 \, ,
\end{equation}   
and the Dirichlet boundary conditions 
\begin{equation}
		\begin{aligned}
			u(x_0 = 0, x_1) &= \frac{1}{2}+ \frac{5}{2} x_1^2 \, ,  \quad u(x_0 = 1, x_1) = 2 + \frac{5}{2} x_1^2 \, , \\
			u(x_0 , x_1 =0) &=  \frac{1}{2}+ \frac{3}{2} x_0^3  \, , \quad  u(x_0, x_1 =1) = 3 + \frac{3}{2} x_0^2  \, .
		\end{aligned}
\end{equation}

\begin{figure}[]
 	\centering
 	\subfloat[][Spectral solution function  \label{fig:2d-poisson-dop}]{
		\includegraphics[width=0.45\textwidth]{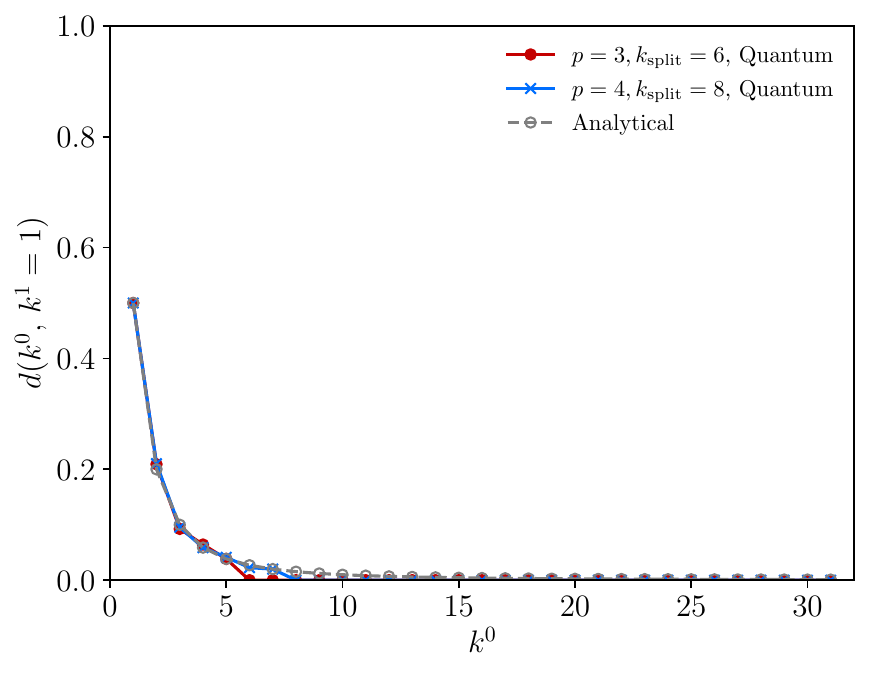} 	}
 	\hfill 
 	\subfloat[][Solution $u(\vec x)$ \label{fig:2d-poisson-soln}]{
 		\includegraphics[width=0.45 \textwidth]{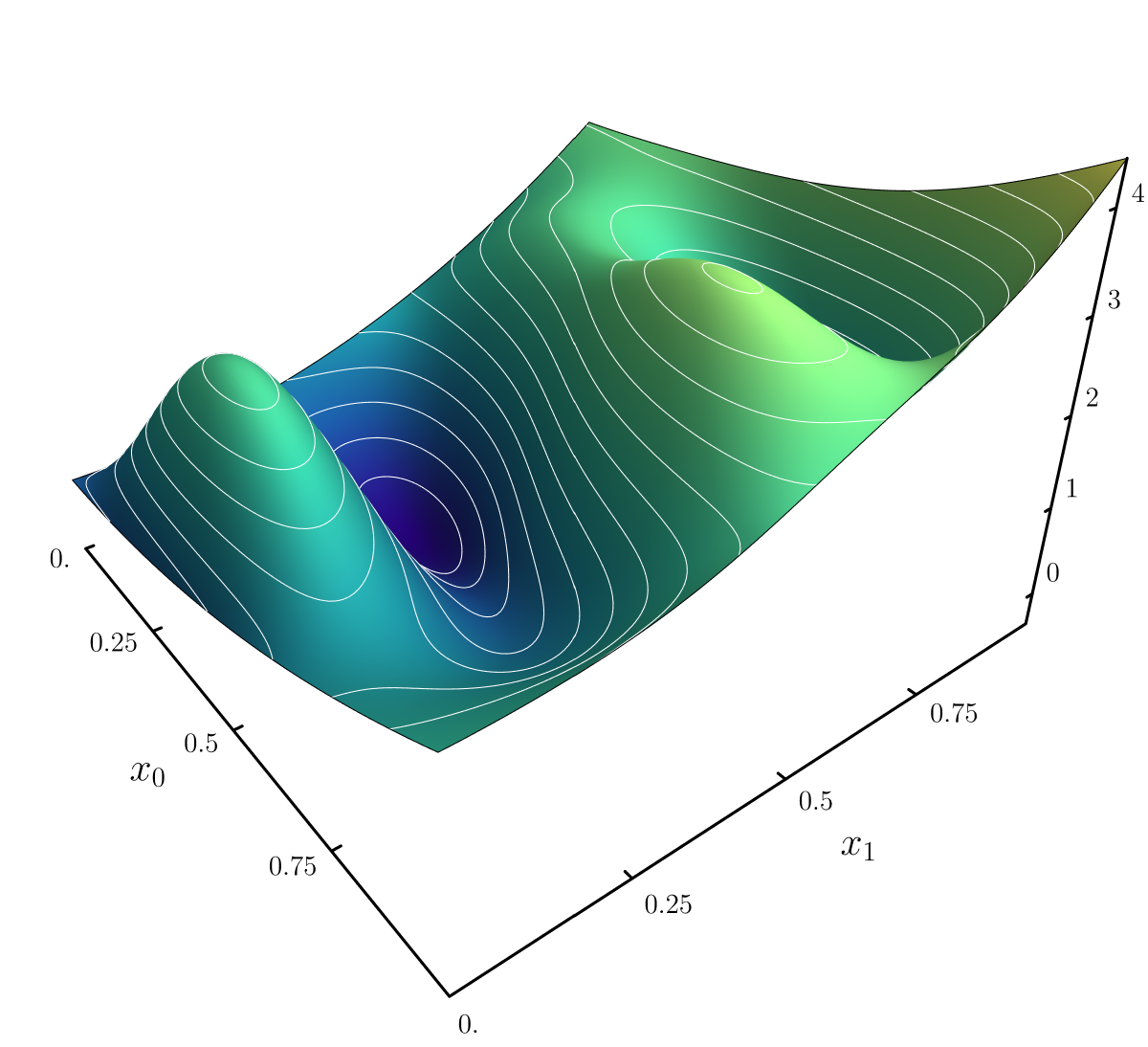}
 	}
 	\caption[]{Two-dimensional inhomogeneous Poisson-Dirichlet problem. In (a), the spectral solution function and its approximation with polynomials of degree~$p=3$ and~$p=4$ are shown.  \label{fig:2d-poisson-nh-contours}}
 \end{figure}

We consider again the decomposition~$u(\vec x) = v (\vec x)+g(\vec x)$, where~$v(\vec x)$ is the solution of the homogeneous Poisson-Dirichlet problem and~$g(\vec x)$ is a sufficiently smooth function which satisfies the boundary conditions. We choose  
\begin{equation}
	g(\vec x) = \frac{1}{2}\left( 3 x_0^2 + 5 x_1^2 + 1 \right) \, .
\end{equation}
In general, this auxiliary function can be, for instance, constructed by approximating the prescribed boundary data using polynomials. The homogeneous Dirichlet-Poisson problem to solve reads 
\begin{equation} \label{eq:homog_dirichlet_inhomog_2d}
\begin{aligned}
	 - \nabla^2 v(\vec x) &= f (\vec x ) + \nabla^2  g(\vec x) \, , \quad  \forall \vec x \in \Omega \, , \\
	 v(\vec x) &= 0  \, , \quad \forall \vec x \in \partial \Omega \,   . 
\end{aligned}
\end{equation}
The bivariate function for computing the Fourier space solution of this problem is given by 
\begin{equation} 
	\label{eq:fourier_sol_dp}
	d(k^0, k^1) = \begin{cases}
		\left ( \dfrac{L}{\pi} \right )^2 \dfrac{1}{(k^0)^2 + (k^1)^2}  \quad & (k^0, k^1) \in (0, N/2)\times (0, N/2) \\
		0 & \text{otherwise} \, . 
	\end{cases}
\end{equation}
As in one-dimensional problems, $d(k^0, k^1)$ is first approximated with a bivariate piecewise polynomial and then encoded as a unitary~$U_P$, see~\ref{appx:bivariate_piecewise}. In two-dimensional problems, the number of polynomials must be restricted to limit the number of gates so that the quantum circuits can still be simulated in Qiskit in a reasonable time. This limitation is relevant to the simulation of the circuit using a classical computer and will not be a problem for future quantum computers. To limit the number of gates, we decompose the Fourier domain $(0, N/2) \times (0, N/2)$ into four subdomains and approximate~$d(k^0, k^1)$ only on the subdomain~$(0, k_\text{split} ) \times (0, k_\text{split})$. In the remaining three subdomains, $d(k^0, k^1)$ is set to zero. We choose multivariate polynomials of degree $p=3$ and $p=4$ as approximants and choose $ k_\text{split} = 6$ for $ p=3$ and $ k _\text{split} = 8$ for $ p=4$.  In Figure~\ref{fig:2d-poisson-dop}, the exact function~$d(k^0, k^1)$ and its two different approximations for~$k^1=1$ are plotted. Notice the very fast decay of the solution function~$d(k^0, k^1=1)$ with increasing~$k^0$, which justifies the neglecting of the higher modes. 

\begin{figure}[]
 	\centering
 	\subfloat[][$L^2$-norm error of $u(x)$ \label{fig:2d-poisson-conv}]{
 		\includegraphics[scale=0.475]{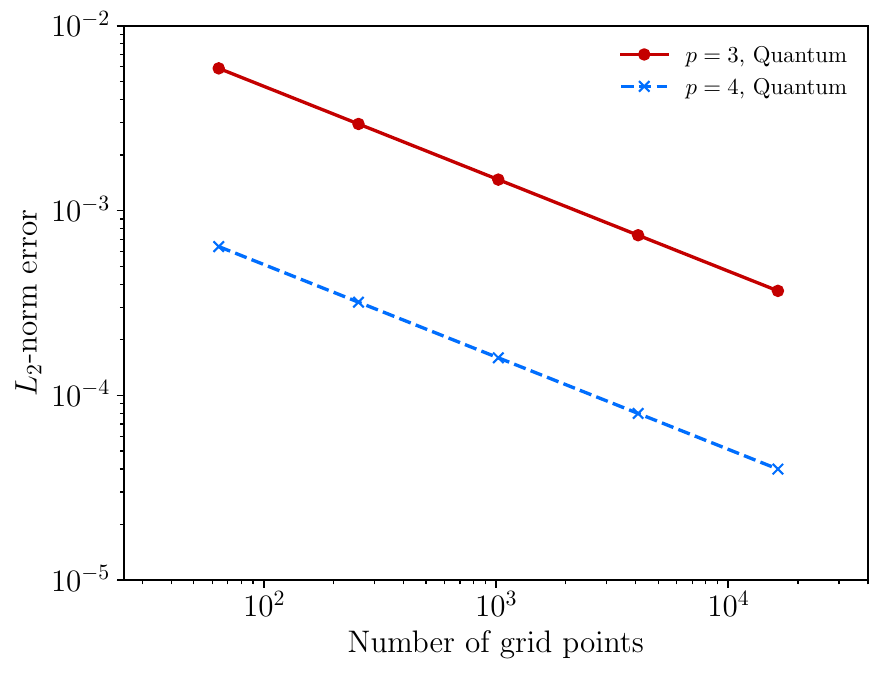}
 	}
 	\hfill
 	\subfloat[][Total number of~$CNOT$ and~$U_3$ gates \label{fig:2d-poisson-cost}]{
 		\includegraphics[scale=0.475]{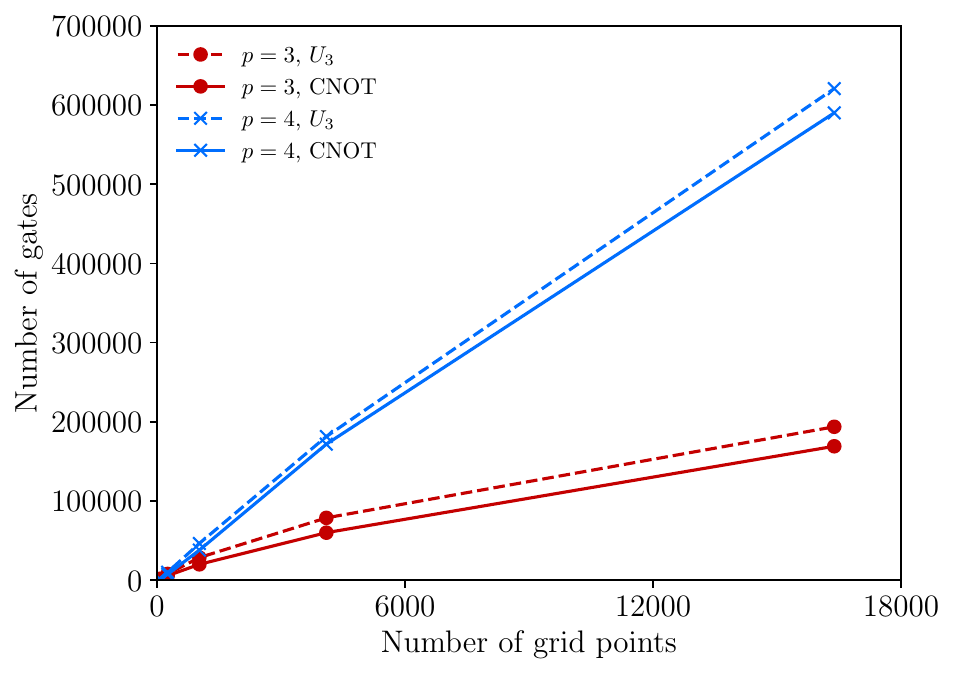}
 	}
 	\caption[]{Two-dimensional inhomogeneous Poisson-Dirichlet problem. \label{fig:2d-poisson-nh}}
 \end{figure}

Figure~\ref{fig:2d-poisson-soln} shows the contour plot of the quantum solution using a discretisation with  $64 \times 64$ grid points and a quartic approximation of the solution function~$d(k^0, k^1)$. The convergence of $L^2$-norm error is depicted in Figure~\ref{fig:2d-poisson-conv}, confirming the uniform convergence of the proposed quantum approach. It is worth noting that the chosen cut-off value~$k_\text{split}$ for the solution function is kept for all discretisations the same. The associated scaling of the total number of~$CNOT$ and~$U_3$ gates is plotted in Figure~\ref{fig:2d-poisson-cost}. Notice the significant increase in the gate counts in comparison to one-dimensional problems. Irrespective of this, the polylogarithmic scaling of the total number of gates is evident. 

\subsubsection{Fractional stochastic partial differential equation} 
\label{sec:twoD-frac}
%
We now apply the proposed quantum approach to the fractional stochastic partial differential representation of two-dimensional random fields. The two-dimensional generalisation of the stochastic ordinary differential equation~\eqref{eq:fractional_ode}  takes the form 
\begin{equation}
		\left ( \kappa^2  - \nabla^2   \right )^\beta  u ( \vec x) = \frac{1}{\tau} f (\vec x) \, ,  \quad \forall  \vec x \in (0, 1) \times  (0, 1) \, , 
\end{equation} 
where~$f(\vec x)$ denotes white noise and~$\kappa$, $\beta$ and~$\tau$ are prescribed parameters. The corresponding Fourier space solution is given by the bivariate function 
\begin{equation} 
	\label{eq:diag_2D_frac}
	d(k^0, k^1) = \begin{cases}
		\dfrac{1}{  \tau  \left  ( \kappa^2 + \dfrac{\pi^2}{L^2}    \left ( (k^0)^2 + (k^1)^2  \right ) \right )^\beta }  \quad & (k^0, k^1) \in (0, N/2)\times (0, N/2) \\
		0 & \text{otherwise} \, . 
	\end{cases}
\end{equation}

\begin{figure}[!h]
	\centering
	\subfloat[][Spectral solution function \label{fig:2d-frac-approx-curve}]{
			\includegraphics[scale=0.475]{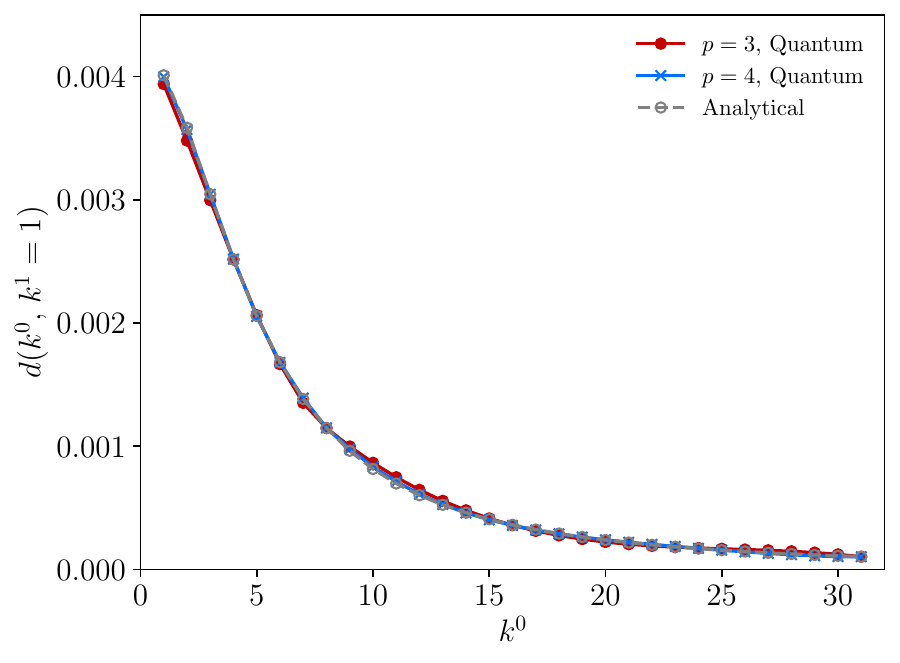}
		}
	\hfill
	\subfloat[][Covariance function  \label{fig:2d-matern-cov}]{
			\includegraphics[scale=0.475]{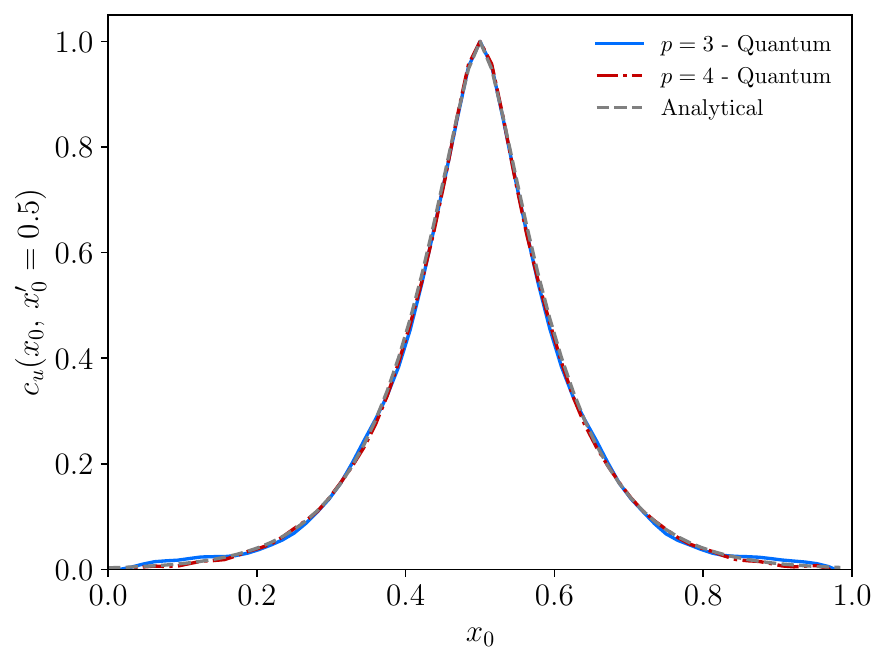}
		}
	\caption[]{Two-dimensional fractional stochastic partial differential equation.  Analytic and approximate spectral solution and Mat\'ern functions for polynomial approximants of degree~$p=3$ and~$p=4$. \label{fig:2d-frac-approx} }
%
	\centering
	\subfloat[][Realisation 1 \label{fig:2d-matern-field-surf}]{
		\includegraphics[scale=0.375]{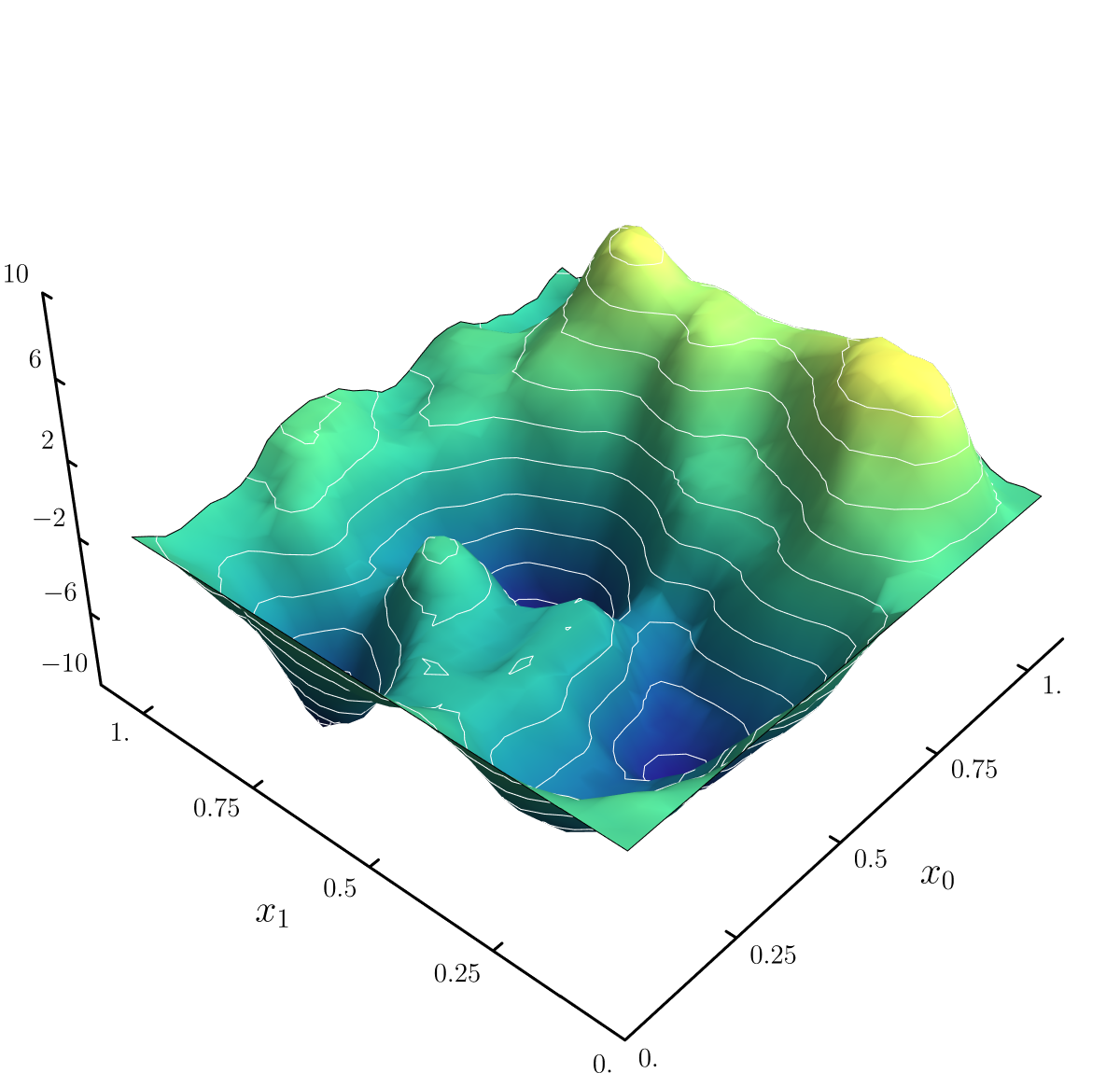}
	}
	\hfill
	\subfloat[][Realisation 2 \label{fig:2d-matern-field-cont}]{
		\includegraphics[scale=0.375]{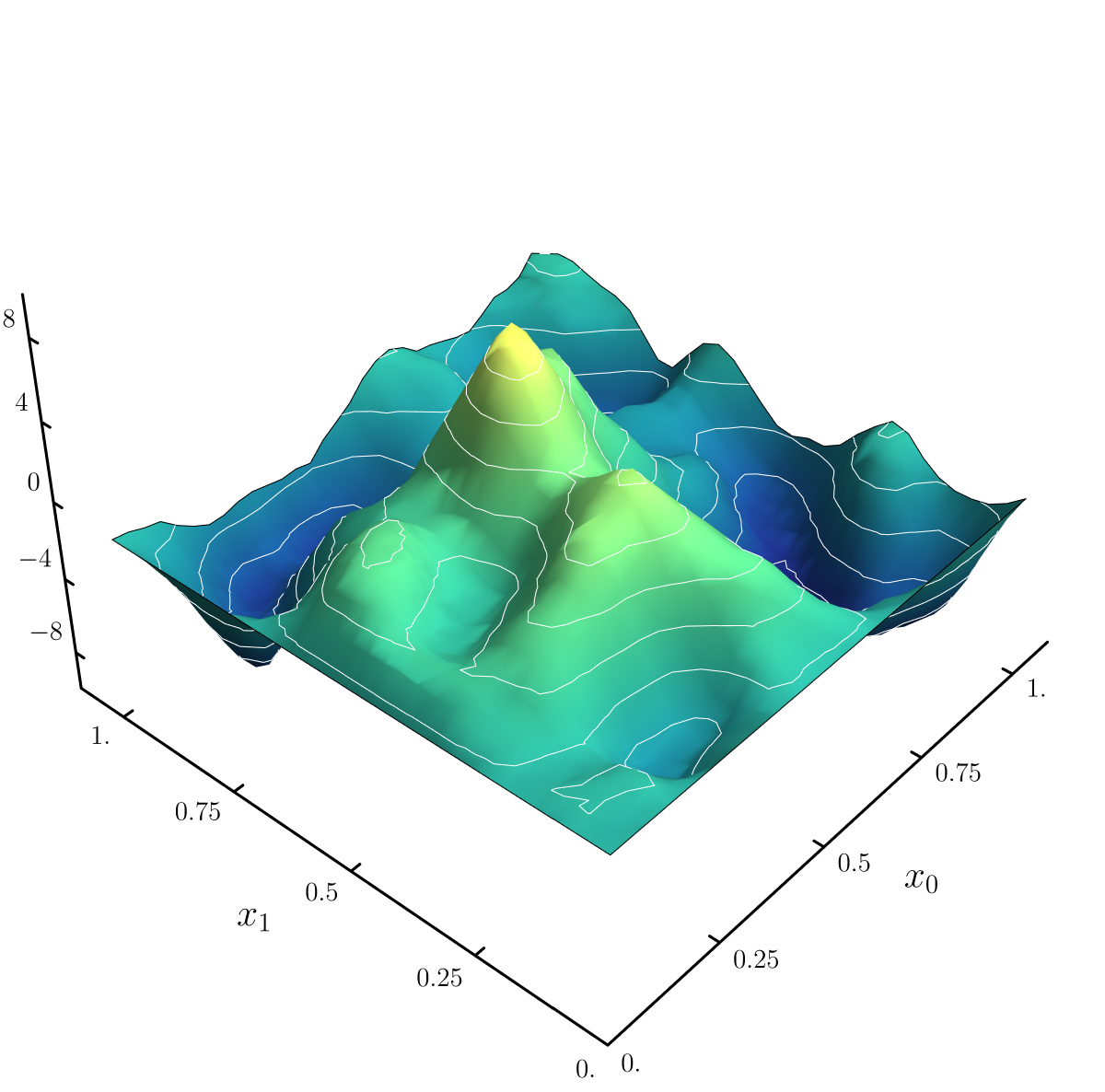}
	}
	\caption[]{Two-dimensional fractional stochastic partial differential equation. Sample solution fields from the Gaussian random field with Mat\'ern covariance. \label{fig:2d-matern-field}}
\end{figure}

The three parameters are chosen as $\kappa=10/\sqrt{2}$, $\tau= 0.01995$, and~$\beta=1$. As shown in Figure~\ref{fig:2d-frac-approx-curve}, the solution function~$d(k^0, k^1)$ is relatively large throughout the Fourier domain. This contrasts with the Fourier space solution function~\eqref{eq:fourier_sol_dp} of the Poisson-Dirichlet problem, which decays rapidly at large wavenumbers. 
We decompose the Fourier domain~$(0,N/2) \times (0,N/2)$ into four subdomains by introducing the wavenumber~$k_\text{split}$ along each axis as discussed in the previous section. Different from the Poisson-Dirichlet problem, we approximate~$d(k^0, k^1)$ in each of the four subdomains with a polynomial. We choose~$k_\text{split}=8$ and bivariate polynomials of degree~$p=3$ and~$p=4$ as approximants. The resulting approximations for $k^1=1$ are compared with the exact~$d(k^0, 1)$ in Figure~\ref{fig:2d-frac-approx-curve}. 

As shown in Figure~\ref{fig:2d-matern-cov} the two different approximations for~$d(k^0, k^1=1)$ yield essentially the same approximation to the covariance function~$c_u(x, x'=1/2)$. The approximate covariance function is given by the entries of the covariance matrix. The good agreement between the quantum approximations and the exact function is evident. Sample realisations of the random field are generated by sampling a white noise forcing vector and solving the partial differential equation. In Figure~\ref{fig:2d-matern-field}, two different realisations of the random field are depicted.

\section{Conclusions}
\label{sec:conclusion}
%
We have introduced a Fourier-based quantum spectral approach to solve non-periodic Dirichlet boundary value problems with polylogarithmic complexity. This represents an exponential speedup relative to the at-best-linear, often worse, complexity achievable with classical approaches.  The quantum speedup is obtained by representing the $N$-dimensional forcing and solution vectors using $n=\log N$ qubits, owing to quantum entanglement, and by designing a unitary solution operator with a polynomial number of one and two-qubit gates. We conjecture that most classical discretisation approaches cannot be efficiently ported to quantum computing due to the inherent difficulty of converting classical vectors and matrices into an exponentially compressed quantum representation. Therefore, amongst the classical methods, well-structured and straightforward discretisation methods, such as the considered Fourier-based spectral approach, are most suited for quantum computing. In the spectral domain, the solution operator simplifies to a diagonal matrix, due to the orthogonality property of the Fourier basis functions, and is relatively easy to invert. Vectors and matrices are converted between the real and spectral domains using the QFT algorithm, which is perhaps the most widely used algorithm throughout quantum computing. As demonstrated, homogeneous Dirichlet boundary conditions can be enforced by adding~$d$ additional qubits to the circuit, where $d$ is the spatial dimension, thereby doubling the domain size and the number of discretisation points. Inhomogeneous Dirichlet boundary conditions are enforced by making use of the superposition principle. 

We close by noting a few possible extensions and algorithmic improvements of the proposed approach. Although we considered only Dirichlet boundary conditions, other boundary conditions can be enforced by choosing different  extensions~\cite{strang1999discrete}. Specifically, Neumann-Dirichlet boundary conditions can be enforced by adding $2d$ qubits and quadrupling the domain size in each coordinate direction. In the present work, our main aim was to demonstrate polylogarithmic complexity irrespective of the total number of gates used. A reduction in gate counts is relevant for near-term quantum computing due to the short coherence times and the limited number of qubits in current quantum computers. The gate counts can be reduced by choosing alternative algorithms for each of the operations in the proposed quantum spectral approach. For instance, it may be possible to replace the used QFT algorithm with an approximate QFT algorithm that omits gates with a very small rotation~\cite{nam2020approximate}. Similarly, it is possible to compute or approximate the inverse of the diagonal solution matrix with alternative techniques, such as QSVT~\cite{martyn2021grand}  or recent Fourier- or Chebyshev-based polynomial approximation techniques \cite{rosenkranz2025quantum}, which may or may not yield a reduction in gate counts.  We did not consider in this paper the extraction of classical data of engineering interest from the obtained quantum states. To this end, recent advances in shadow tomography~\cite{huang2020predicting} and quantum machine learning~\cite{williams2024addressing} may provide practical solutions and warrant further investigation.  In terms of future applications, solving non-periodic boundary value problems on $d$-dimensional hyperrectangles is, among others, relevant in homogenisation~\cite{morin2024fast,risthaus2024fft,paux2025discrete}. Consequently, our quantum solver can be used as a component in the quantum homogenisation approach introduced earlier by Liu et al.~\cite{liu2024towards}. Finally, quantum approaches offer distinct advantages when applied to stochastic problems that require the repeated solution of the same problem with varying inputs~\cite{ deiml2025quantum}. The inputs are usually samples from a random field, which, as demonstrated in this paper, can be obtained by solving a stochastic differential equation. 

\section*{Acknowledgments}
EF is supported through Reward for Excellence Award of the Early Career Development Program at the College of Science and Engineering, University of Glasgow.

\appendix
\section{Algorithmic complexity \label{appx:alg_complex}}
%
We assess the complexity of quantum algorithms by expressing circuits using the gate set~$\{ CNOT, \, U_3\}$ with the two-qubit gate
\begin{equation}
	CNOT = \ket 0 \bra 0 \otimes I + \ket 1 \bra 1 \otimes X 
	= 
	\begin{pmatrix}
		1 & 0 & 0 & 0 \\ 
		0 & 1 & 0 & 0 \\ 
		0 & 0 & 0 & 1 \\ 
		0 & 0 & 1 & 0 
	\end{pmatrix} \, ,
\end{equation}
and the single qubit gate
\begin{equation}
	U_3 (\theta, \phi, \lambda) = \begin{pmatrix}  \cos \frac{\theta}{2} & - e^{i \lambda} \sin \frac{\theta}{2} \\[0.2em] e^{i \phi} \sin \frac{\theta}{2}  & e^{i(\phi+\lambda)} \cos \frac{\theta}{2} \end{pmatrix} \, ,
\end{equation}
where~$\theta, \phi, \lambda \in \mathbb R$ are three angles. The gate set~$\{ CNOT, \, U_3\}$ is universal in the sense that any quantum circuit can be expressed using its two gates. 

All of our circuits are implemented in Qiskit~\cite{qiskit2024}. Subsequently, they are expressed using only~$\{ CNOT, \, U_3\}$ gates by compiling them either in Qiskit or TKET~\cite{tket2020}. In Qiskit, we use the \texttt{transpile} command, and in TKET, the~\texttt{qiskit\_to\_tk} and \texttt{get\_compiled\_circuit} commands, with optimisation level~3. Usually, the number of gates in TKET is significantly lower than in Qiskit. 

\section{Encoding of bivariate piecewise polynomials \label{appx:bivariate_piecewise}}
%
Consider the two-dimensional domain~$\Omega \subset \mathbb N^2$ and its partitioning into~$M\times M$ non-overlapping subdomains~$\omega_{\vec c}$ such that 
\begin{equation}
	\Omega = \bigcup_{\vec c} \omega_{\vec c} \, , \quad \vec c = (c^0, c^1), \, \quad c^0,c^1  = 0, \, \dotsc, \, M-1 \, .
\end{equation}
Each subdomain~$\omega_{\vec c}$ is defined in terms of its bounding box 
\begin{equation}
	\omega_{\vec c} = \left [ l^0_{\vec c}, u_{\vec c}^0 \right ) \times \left [ l^1_{\vec c}, u_{\vec c}^1 \right ) \, ,
\end{equation}
where~$l^0_{\vec c}$ and~$u^0_{\vec c}$ are the lower and upper grid points in the coordinate direction~$k^0$, and~$l^1_{\vec c}$ and~$u^1_{\vec c}$ are the respective grid points in the coordinate direction~$k^1$. Each subdomain~$\omega_{\vec c}$ is associated with a bivariate polynomial
\begin{equation}
	\tilde d_{\vec c} (k^0, k^1)= \begin{pmatrix} 1 \\ k^0 \\ (k^0)^2 \\ \vdots  \\  (k^0)^p  \end{pmatrix}^\trans   
	 \begin{pmatrix} \alpha_{\vec c, 00} & \alpha_{\vec c, 01} & \alpha_{\vec c, 02} & \cdots & \alpha_{\vec c, 0p}  \\
	 \alpha_{\vec c, 10} & \alpha_{\vec c, 11} & \alpha_{\vec c, 12} & \dotsc & \alpha_{\vec c, 1p} \\
	  \alpha_{\vec c, 20} & \alpha_{\vec c, 21} & \alpha_{\vec c, 22} & \dotsc & \alpha_{\vec c, 2p} \\
	   \vdots & \vdots & \vdots & \vdots & \vdots \\
	   \alpha_{\vec c, p0} & \alpha_{\vec c, p1} & \alpha_{\vec c, p2} & \dotsc & \alpha_{\vec c, pp}
	   \end{pmatrix}  \begin{pmatrix} 1 \\ k^1 \\ (k^1)^2 \\ \vdots \\ (k^1)^p  \end{pmatrix}  
	 =   \begin{pmatrix} 1 \\ k^0 \\ (k^0)^2 \\ \vdots  \\  (k^0)^p  \end{pmatrix}^\trans   \vec \alpha_{\vec c} \begin{pmatrix} 1 \\ k^1 \\ (k^1)^2 \\ \vdots \\ (k^1)^p  \end{pmatrix} \, ,
\end{equation}
where~$p$ is the polynomial degree and $\vec \alpha_{\vec c} \in \mathbb R^{(p+1) \times (p+1)}$ is the coefficient matrix.

The support of~$\tilde d_{\vec c} (k^0, k^1)$ is restricted to the subdomain~$\omega_{\vec c}$ with the same index~$\vec c$. Hence, to evaluate it at a given grid point~$\vec k=(k^0, k^1)$ first the subdomain index~$\vec c$ must be found, such that \mbox{$ l^0_{\vec c} \le k^0 <  u^0_{\vec c}$}  and \mbox{$ l^1_{\vec c} \le k^1 <  u^1_{\vec c} $}. Because cell boundaries are axis-aligned, the subdomain index~$\vec c$ can be found using~$4M$ integer comparisons.  The number of integer comparisons can be reduced by a factor of two by only identifying the set of all subdomain indices~$\{ \vec c\}$ such that~$k^0<u^0_{\vec c}$ and~$k^1<u^1_{\vec c}$ and evaluating all the associated polynomials~$\{ d_{\vec c}(\vec k)\}$. To facilitate this, the polynomials are modified as follows
\begin{equation} \label{eq:bivar_mod_coeffs}
\begin{aligned}
\tilde d_{(0,0)} (\vec k) & \leftarrow \tilde d_{(0,0)} (\vec k)  \\ 
\tilde d_{(1,0)} (\vec k) & \leftarrow \tilde d_{(1,0)} (\vec k) - \tilde d_{(0,0)}  (\vec k)  \\ 
\tilde d_{(2,0)} (\vec k) & \leftarrow \tilde d_{(2,0)} (\vec k) - \tilde d_{(1,0)} (\vec k)   \\ 
\vdots  \\
\tilde d_{(0,1)} (\vec k) & \leftarrow \tilde d_{(0,1)} (\vec k) - \tilde d_{(0,0)}  (\vec k) \\ 
\tilde d_{(1,1)} (\vec k) & \leftarrow  \tilde d_{(1,1)} (\vec k) - \tilde d_{(0,1)} (\vec k)  - \tilde d_{(1,0)} (\vec k) +  \tilde d_{(0,0)}  (\vec k) \\ 
\tilde d_{(2,1)} (\vec k) & \leftarrow  \tilde d_{(2,1)} (\vec k) - \tilde d_{(1,1)} (\vec k)  - \tilde d_{(2,0)} (\vec k) +  \tilde d_{(1,0)}  (\vec k) \\ 
\vdots \, . 
\end{aligned}
\end{equation}

\begin{figure}
	\centering
	\scalebox{0.85}{
		\input{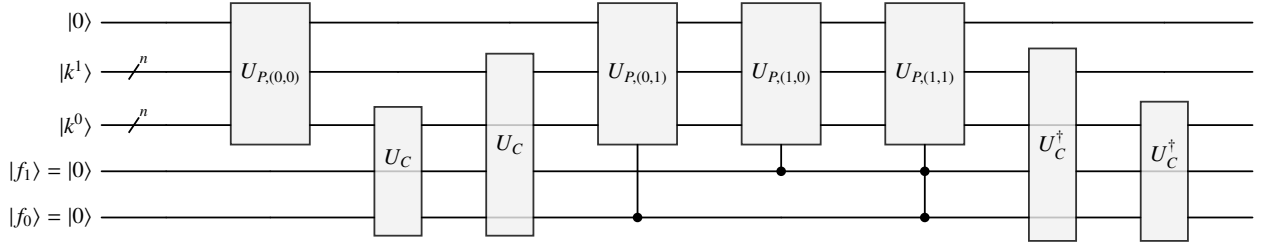}
	}
	\caption{Quantum circuit for encoding a piecewise bivariate polynomial. The domain~$[0, 2^n-1]^2$ is partitioned into $2\times 2$ subdomains of equal size. The values of the two flag qubits, $\ket {f_0}$ and $\ket {f_1}$, are set by the integer comparator~$U_C$ and determine the relevant polynomials for a given point evaluation point~$\vec k= (k^0, k^1) $. The first flag is in state $\ket {f_0} = \ket{1}$ when $k^0 \ge 2{n-1}$ and the second flag is in state~$\ket{f_1}=\ket 1$ when $k^0 \ge 2{n-1}$. After the circuit is executed the function value is the amplitude of the state~$\ket 1 \ket{k^1} \ket{k^0} \ket 0 \ket 0$. \label{fig:piecewise_encoding}}
\end{figure}

As an example, the domain~$\Omega =  [ 0, 2^n-1]^2$ consists of four subdomains~$\omega_{\vec c}$ so that $\vec c \in \{ (0,0), \, (0,1), \, (1,0), \, (1,1) \}$.  The four non-overlapping subdomains are of equal size so that, e.g., the subdomain~$\omega_{(0,0)}$  has the domain $[0, 2^{n-1}-1)^2$. The piecewise polynomial defined over~$\Omega$ is evaluated with the circuit shown in Figure~\ref{fig:piecewise_encoding}.  The circuit consists of four unitaries~$U_{P, \vec c}$  for evaluating the four polynomials defined according to~\eqref{eq:bivar_mod_coeffs}. The unitaries are conditioned on the two ancilla qubits~$f_0$ and~$f_1$, which are set by the integer comparator~$U_C(2^{n-1})$ as applied to the first and second coordinate direction~$k^0$ and~$k^1$, i.e.,   
\begin{equation}
	U_C(2^{n-1}) \colon  \, \ket {k^0}  \ket {f_0 = 0}  \mapsto  \ket {k^0}  \ket{k^0\ge 2^{n-1}}  \, \quad  U_C(2^{n-1})  \colon  \, \ket {k^1}  \ket {f_1 = 0}  \mapsto  \ket {k^1}  \ket{k^1 \ge 2^{n-1}}  \, .
\end{equation}
Focusing on the first comparator unitary, it has the output~$\ket {k^0} \ket 1$ for all grid points with~$k^0 \ge 2^{n-1}$ and~$\ket {k^0} \ket 0$ otherwise.  Most quantum computing libraries, including Qiskit, provide an implementation of $ U_C$ with polylogarithmic complexity. $U_C$ can be, for instance, defined in terms of binary addition.

\section{Stochastic differential equation representation of random fields \label{appx:ode_matern}}
%
A Mat\'ern-type random field~$u(x)$ is the solution of the stochastic ordinary differential equation 
\begin{equation}  \label{eq:appx_fractional_ode_ex}
		\left ( \kappa^2  - \frac{\D^2 }{\D x^2} \right )^\beta  u ( x) = \frac{1}{\tau} f (x) \, ,  \quad  \forall x \in (-\infty, +\infty) \, ,
\end{equation} 
where the Gaussian white-noise forcing~$f(x)$ is given by 
\begin{equation}
	f(x) \sim \mathcal{GP}\left( 0, \delta (x - x') \right) \, , 
\end{equation}
and has the covariance function
\begin{equation}
	\cov ( f(x), f(x'))=  c_{f} (x,x') = \expect [f(x) f(x')] = \delta(x-x')   \, .
\end{equation}
The  solution of~\eqref{eq:appx_fractional_ode_ex} for a single realisation of~$f(x)$ can be determined via the Greens function~$g( x, x')$ as 
\begin{equation}
	u(x) = \int g(x, x') f(x') \D x \, .
\end{equation}
By linearity of expectation, the random solution field~$u(x)$ has the covariance
\begin{equation}
	\cov(u(x),u(x')) = \expect[u(x)u(x')] = \int \int g(x,x'') \expect [f(x'') f(x''')] g(x', x''')  \D x'' \D x'''  
	=   \int g(x,x'') g(x', x'') \D x'' \, .
\end{equation}
This expression is equal to the classical Mat\'ern covariance  given by 
\begin{equation} \label{eq:matern-cov}
	 \cov \left(u(x) , u(x') \right) = c_u(x, x')  = \frac{1}{2^{\nu - 1} \, \Gamma(\nu)} \left( \frac{\sqrt{2 \nu}}{\ell} \| x - x'\|\right)^{\nu} K_{\nu} \left( \frac{\sqrt{2 \nu}}{\ell} \| x - x'\|\right) \, .
\end{equation}
The parameters $\kappa, \, \beta$, and $\tau$ of the differential equation are related to the Mat\'ern covariance  parameters as follows
\begin{equation}
	\label{eq:spde-params}
	\kappa = \frac{\sqrt{2 \nu}}{\ell} \, , \quad \beta = \frac{\nu}{2} +  \frac{d}{4} \, , \quad \tau^2 = \frac{\Gamma (\nu)}{\Gamma(\nu + d/2) (4 \pi)^{d/2} \kappa^{2 \nu}} \, ,
\end{equation}
where $\Gamma(\cdot)$ is the Gamma function, and $K_{\nu}(\cdot)$ is the modified Bessel function of the second kind and $d$ denotes the spatial dimension.

Finally, since the stochastic differential equation is linear  and the random forcing is Gaussian, the solution is also a Gaussian field, which can be compactly written as 
\begin{equation} \label{eq:u_random_field}
	u(x) \sim \mathcal{GP}\left( 0, \frac{1}{\tau^2}  c_u (x, x') \right) \, .
\end{equation}
For further details see~\cite{koh2023stochastic,lindgren2011explicit}.

\bibliographystyle{elsarticle-num-names}
\bibliography{quantumEllipticPDE}

\end{document}